\documentclass[a4paper, DIV10,12pt]{article}
\usepackage[a4paper, left=2.5cm, right=2.5cm, top=2.5cm, bottom=2cm]{geometry}
\usepackage[utf8]{inputenc}
\usepackage{lmodern} 
\usepackage{todonotes}
\usepackage[english]{babel}
\usepackage{amsmath}
\usepackage{algorithmic}
\usepackage[T1]{fontenc}
\numberwithin{equation}{section}
\usepackage{amsfonts}
\usepackage{amssymb}
\usepackage[ruled,vlined]{algorithm2e}
\usepackage{graphicx}

\newtheorem{remark}{Remark}[section]

\title{Multilevel Sequential Importance Sampling for Rare Event Estimation}
\author{F. Wagner, J. Latz, I. Papaioannou, E. Ullmann}
\date{\today}

\begin{document}

\maketitle

\begin{abstract}
The estimation of the probability of rare events is an important task in reliability and risk assessment. We consider failure events that are expressed in terms of a limit state function, which depends on the solution of a partial differential equation (PDE). Since numerical evaluations of PDEs are computationally expensive, estimating such probabilities of failure by Monte Carlo sampling is intractable. We develop a novel estimator based on a Sequential Importance sampler using discretizations of PDE-based limit state functions with different accuracies. A twofold adaptive algorithm ensures that we obtain an estimate based on the desired discretization accuracy. Moreover, we suggest and study the choice of the MCMC kernel for use with Sequential Importance sampling. Instead of the popular adaptive conditional sampling method, we propose a new algorithm that uses independent proposals from an adaptively constructed von Mises-Fisher-Nakagami distribution.
\vspace{0.5cm}
\\\textbf{Keywords:} Reliability analysis, Importance Sampling, Multilevel Monte Carlo, Subset Simulation, Markov Chain Monte Carlo
\end{abstract}

\section{Introduction}\label{Sec: Introduction}
Estimating the probability of rare events is crucial in reliability analysis and risk management and arises in applications in many fields. For instance, the authors in \cite{Agarwal17} examine rare events arising in financial risk settings while \cite{Morio15} studies the probability of collision between space debris and satellites. In planning a radioactive waste repository \cite{Cornaton08,Noseck08}, one is interested in the probability that radioactive particles leave the repository and pollute the groundwater in a long time horizon. The particle flow can be simulated by a finite element (FEM) \cite{Braess07} approximation of the groundwater flow and transport equation. Since the subsurface properties of the whole domain of interest are uncertain or only measurable at finitely many points, the soil is modelled as a random field. The particle transport has to be simulated for various realizations of the random field to estimate the probability that the radioactive particles come back to the human environment, which is a rare event.
\\All applications have in common that the probabilities of the events are small $(<10^{-4})$ and the \emph{limit state function} (LSF) underlies a computationally demanding model which depends on the discretization of the domain. If the discretization level is high, i.e., the mesh size is small, the FEM approximation is accurate but also cost intensive. These issues complicate the estimation of the probability of failure.
\\Before we introduce our novel approach, we give a brief overview of existing algorithms. On the one hand, there are deterministic approximation methods, such as the \emph{first} and \emph{second order reliability method} (FORM, SORM) \cite{Melchers18}, which aim at approximating the domain of parameters which lead to failure events. On the other hand, there are sampling based methods, which approximate the probability of failure events. Unlike approximation methods, sampling approaches are based on sample estimates and are usually more robust in terms of the complexity of the LSF. Since our novel approach is a sampling method, we focus on this category and give a larger overview. 
\\\emph{Monte Carlo sampling} \cite{Fishman96,Rubinstein16} can be easily applied to estimate the probability of failure and yields an unbiased estimator. However, due to the mentioned issues of rare event settings, the Monte Carlo estimator becomes intractable, since hardly any sample contributes to the rare event and each sample requires a cost intensive function evaluation. Therefore, variance reduction techniques have been developed to reduce the number of samples for obtaining an accurate estimate. For instance, the idea of \emph{Multilevel Splitting} \cite{Botev12,Glasserman99} and \emph{Subset Simulation} (SuS) \cite{Au01,Au14} is to decompose the rare event into a sequence of nested events. This enables expressing the probability of the rare event as a product of conditional probabilities of more frequent events. These methods require sampling from a sequence of probability density functions which is achieved with \emph{Markov chain Monte Carlo} (MCMC) methods \cite{Papaioannou15,Wang19}.
\\\emph{Importance sampling} (IS) methods employ an alternative sampling density, which if chosen properly can reduce considerably the variance of the standard Monte Carlo estimator \cite{Kahn53}. The optimal choice of the sampling density is the density of the input variables conditional on the failure domain. However, direct sampling from the optimal density is not feasible, because the location of the failure domain is unknown prior to performing the simulation. As in Multilevel Splitting or SuS, a sequential approach can be applied to approximate the optimal IS density in a sequential manner. This leads to \emph{Sequential Importance Sampling} (SIS) \cite{Papaioannou18,Papaioannou16} or \emph{Sequential Monte Carlo} (SMC) \cite{Cerou12} for the estimation of rare events. In our novel approach, we consider an adaptive methodology similar to adaptive SMC \cite{Beskos13,Doucet09,Jasra11}. Another approach to estimate the optimal sampling density sequentially is the \emph{Cross-entropy} method \cite{Geyer19}, where the sampling density minimizes the Kullback-Leibler divergence to the optimal density within a family of parametrized densities. IS can also be applied to a hyperplane that is perpendicular to an important direction, a method known as \emph{line sampling} \cite{Angelis15,Koutsourelakis04,Rackwitz01}.
\\The previous algorithms have the drawback that all evaluations have to be performed with respect to the same LSF. The evaluation of the LSF could require the solution of a discretized PDE, which depends on the mesh size of the computational domain. Since computational costs increase with decreasing mesh size, we wish to construct a method wherein the discretized PDE is solved on fine meshes only for very few realizations. Therefore, we apply a multilevel approach that uses a hierarchy of discretization levels. The authors in \cite{Elfverson16} use the telescoping sum approach of \cite{Giles15} to estimate the probability of failure. Applying the multilevel idea to the previously described methods gives \emph{Multilevel Subset Simulation} (MLSuS) \cite{Ullmann15} and \emph{Multilevel Sequential Monte Carlo} \cite{Beskos17,Moral17}. Moreover, a multi-fidelity approach combined with the cross-entropy method is investigated in \cite{Peherstorfer18}. Furthermore, the work in \cite{Latz18} develops the \emph{Multilevel Sequential$^2$ Monte Carlo} ($\mathrm{MLS}^2\mathrm{MC}$) estimator, which is a twofold sequential algorithm for Bayesian inverse problems. 
\\In this paper, we consider SuS and SIS as well as their multilevel versions. In more detail, an MCMC algorithm \cite{Cotter13,Hastings70} is applied within SuS to gradually shift samples into consecutive domains, which are defined by the sequence of nested events. By the \emph{nestedness property} \cite{Papaioannou15}, the simulated Markov chains do not require a burn-in period, since seeds are already distributed approximately according to the target distribution. Therefore, SuS is an efficient but slightly biased estimator \cite{Au01}. The MLSuS method, given in \cite{Ullmann15}, employs a hierarchy of discretization levels and enables the usage of coarse grid function evaluations. MLSuS saves significant computational costs compared to SuS if the failure domains between discretization levels are still nested. However, nestedness is no longer guaranteed in the multilevel setting since the sequence of consecutive domains is based on LSFs with different accuracies. Therefore, a second MCMC step has to be performed. Additionally, a burn-in period is proposed since seeds are no longer distributed (approximately) according to the target distribution. Both issues increase the computational costs of the MLSuS estimator; and thus decrease its efficiency. However, a level dependent parameter dimension can be applied to reduce variances between two accuracy levels of the LSF and approximately satisfy the nestedness property.
\\The nestedness issue of MLSuS is our main motivation to implement the $\mathrm{MLS}^2\mathrm{MC}$ algorithm for rare event estimation. Nestedness is not an issue for $\mathrm{MLS}^2\mathrm{MC}$; the method samples a sequence of non-zero densities with IS and chooses each IS density to be close to each target density in the sequence. The idea of the $\mathrm{MLS}^2\mathrm{MC}$ method is combined with the SIS approach and yields a \emph{Multilevel Sequential Importance Sampling} (MLSIS) estimator for rare events. Note that both MLSIS as well as MLSuS are not based on the telescoping sum approach. To achieve an even more efficient algorithm, we apply the level dependent parameter dimension approach of \cite{Ullmann15}. As SIS, the MLSIS method requires an MCMC algorithm to shift samples into consecutive target distributions. We consider an independent MCMC sampler that uses the \emph{von Mises-Fisher Nakagami} (vMFN) distribution model fitted with the available weighted samples at each sampling level as the proposal distribution. The vMFN distribution is applied in \cite{Papaioannou19} as a parametrized family of probability distributions for the Cross-entropy method which yields an efficient algorithm even in high dimensions. Employing the vMFN distribution as a proposal density is another main contribution of our work.
\\The paper is structured as follows. In Section \ref{Section SIS}, the problem setting of estimating the probability of failure is defined and SIS as well as SuS are explained. The MLSIS estimator is described in Section \ref{Section MLSIS}. In Section \ref{Section MCMC}, two MCMC algorithms are studied which are applied within SIS and MLSIS. In Section \ref{chapter numerical experiments}, the studied estimators are applied to 1D and 2D test problems and the MLSIS estimator is compared with SIS as well as SuS and MLSuS. In Section \ref{Section Conclusion}, a summary of the discussion and an outlook are given.

\section{Background}\label{Section SIS}

\subsection{Problem Setting}
Consider the probability space $(\Omega, \mathcal{F}, \mathbb{P})$ and a random variable $U:\Omega\rightarrow\mathbb{R}^n$. By \cite{Kiureghian86,Hohenbichler81} it is assumed, without loss of generality, that $U$ is distributed according to the $n$-variate standard normal distribution with density function $\varphi_n$. If a non-Gaussian random variable $\tilde{U}$ is used, an isoprobabilistic transformation $U=T(\tilde{U})$ is applied. Failure is defined in terms of an LSF $G:\mathbb{R}^n\rightarrow\mathbb{R}$ such that $G(U(\omega))\le 0$ for $\omega\in\Omega$. In many applications, the LSF $G$ is not analytically given. We can only evaluate an approximation $G_{\ell}$, where $\ell$ represents the discretization level. Increasing $\ell$ leads to a more accurate approximation. In the numerical examples presented in this paper, $G_{\ell}$ requires the solution of a PDE and $\ell$ specifies the mesh size of an FEM approximation. The \emph{probability of failure} is defined as the measure of the \emph{failure domain} $A:= \{\omega\in\Omega: G(U(\omega))\le 0\}$, which is expressed as
\begin{align}
P_{f} := \mathbb{P}[A] = \mathbb{P}[G(U)\le0] = \int_{G(u) \le 0} \varphi_n(u) \mathrm{d}u.\label{probability of failure}
\end{align}
Using $G_{\ell}$ instead of $G$ in (\ref{probability of failure}) gives the approximation $P_{f,\ell}$, which includes numerical errors due to approximating the exact LSF $G$. Convergence is expected for increasing the level $\ell$, i.e., decreasing the finite element mesh size.
\\ The probability of failure can be estimated by crude Monte Carlo sampling \cite{Fishman96}. By evaluating $G_{\ell}$ on the discretization level $\ell$ for $N\in\mathbb{N}$ independent samples distributed according to $\varphi_n$, we obtain the (single-level) Monte Carlo estimator $\hat{P}_{f,\ell}^{\mathrm{MC}}$ for $P_{f,\ell}$
\begin{align}
\hat{P}_{f,\ell}^{\mathrm{MC}} = \frac{1}{N} \sum_{k=1}^N I\left(G_{\ell}(u_k)\le 0\right),\label{MC esimator}
\end{align}
where $I$ denotes the indicator function; i.e., $I(\mathrm{true}) = 1$ and $I(\mathrm{false}) = 0$. $\hat{P}_{f,\ell}^{\mathrm{MC}}$ is an unbiased estimator and easy to implement. Since the coefficient of variation of $\hat{P}_{f,\ell}^{\mathrm{MC}}$ is inversely proportional to the probability of failure $P_{f,\ell}$, see \cite{Papaioannou16}, a large number of samples is required if $P_{f,\ell}$ is small and a small coefficient of variation should be achieved. Hence, huge computational costs are required if $G_{\ell}$ is a cost demanding evaluation. This makes crude Monte Carlo sampling impractical for the estimation of rare failure probabilities. 

\subsection{Subset Simulation and Multilevel Subset Simulation}\label{SuS and MLSuS}
SuS and MLSuS are alternative approaches where the failure probability is estimated by a product of conditional probabilities. Consider the sequence of domains $B_0, B_1,\dots,B_S$, where $B_S = A$ is the failure domain. In both approaches, the sequence is constructed such that
\begin{align}
P(B_j\mid B_{j-1})= \hat{p}_0 \in(0,1),\label{SuS parameter}
\end{align}
while $\hat{p}_0$ is chosen to ensure that samples of $B_j$ can be easily generated from samples of $B_{j-1}$ \cite{Au01}. In SuS, the sequence of domains is nested, i.e., $B_{j}\subset B_{j-1}$ for $j=1,\dots,S$, since the discretization level is fixed. Hence, the SuS estimator is given as 
\begin{align*}
\hat{P}_{f,\ell}^{\mathrm{SuS}} := \hat{P}_{B_1} \prod_{j=2}^S \hat{P}_{B_j\mid B_{j-1}},
\end{align*}
where $\hat{P}_{B_j\mid B_{j-1}}$ is an estimator for $P(B_j\mid B_{j-1})$. It has been shown in \cite{Papaioannou16} that SuS is a special case of SIS, where the IS densities $p_{j,\ell}$ are chosen as the optimal IS density with respect to the domain $B_j$. In MLSuS \cite{Ullmann15}, the sequence of domains is no longer nested since the domains $B_j$ are defined for different LSFs $G_{\ell}$, in case of a level update. To overcome this problem, the conditional probability $P(B_{j-1}\mid B_{j})$ has to be estimated. This leads to the MLSuS estimator
\begin{align}
\hat{P}_{f,\ell}^{\mathrm{MLSuS}} := \hat{P}_{B_1} \prod_{j=2}^S \frac{\hat{P}_{B_{j}\mid B_{j-1}}}{\hat{P}_{B_{j-1}\mid B_{j}}}.\label{MLSuS estimator}
\end{align}
Moreover, samples which are taken as seeds in the MCMC step are not distributed according to the target distribution. Therefore, a burn-in is required. Both issues lead to increasing computational costs. Note that if the domains $B_j$ for $j=1,..,S$ were nested, then the denominator in (\ref{MLSuS estimator}) is equal to one and no estimator for the denominator is required. To increase the denominator in (\ref{MLSuS estimator}), the authors in \cite{Ullmann15} apply a level dependent parameter dimension. This reduces the variance of two consecutive levels and makes the MLSuS algorithm more robust. In Section \ref{level dependent dim} of this work, we consider a level dependent parameter dimension for the MLSIS algorithm to reduce the variance of consecutive levels. 

\subsection{Importance Sampling}
IS is a variance reduction technique \cite{Agapiou17,Rubinstein16}, where the integral in (\ref{probability of failure}) is calculated with respect to a certain \emph{IS density} $p_{\ell}$. If $p_{\ell}$ takes on large values in the failure domain, many samples following $p_{\ell}$ represent failure events. Therefore, less samples are required to estimate the probability of failure accurately. By \cite{Papaioannou18} the failure probability $P_{f,\ell}$ is expressed as
\begin{align*}
P_{f,\ell} = \int_{\mathbb{R}^n} I\left(G_{\ell}(u)\le 0\right) w_{\ell}(u)p_{\ell}(u)\mathrm{d}u = \mathbb{E}_{p_{\ell}}[I\left(G_{\ell}(u)\le 0\right)w_{\ell}(u)],
\end{align*}
where the \emph{importance weight} is defined as $w_{\ell}(u) := {\varphi_n(u)}/{p_{\ell}(u)}$. Again, crude Monte Carlo sampling is applied, which yields the estimator
\begin{align*}
\hat{P}_{f,\ell}^{\mathrm{IS}} = \frac{1}{N} \sum_{k=1}^{N} I\left(G_{\ell}(u_k)\le 0\right)w_{\ell}(u_k),
\end{align*}
where the samples $\{u_k\}_{k=1}^{N}$ are distributed according to the IS density $p_{\ell}$. $\hat{P}_{f,\ell}^{\mathrm{IS}}$ is an unbiased estimator for $P_{f,\ell}$ if the support of $p_{\ell}$ contains the failure domain $A_{\ell}:= \{\omega\in\Omega: G_{\ell}(U(\omega))\le 0\}$. The \emph{optimal IS density} is given by
\begin{align}
p_{\mathrm{opt},\ell}(u) := \frac{1}{P_{f,\ell}}I(G_{\ell}(u)\le 0) \varphi_n(u),\label{optimal IS density}
\end{align}
which leads to a zero-variance estimator. Since $P_{f,\ell}$ and $A_{\ell}$ are unknown, $p_{\mathrm{opt},\ell}$ cannot be used in practice. In contrast, SIS achieves an approximation of $p_{\mathrm{opt},\ell}$ by approximating the optimal IS distribution in a sequential manner while starting from a known \emph{prior density} $p_0$. 
 
\subsection{Sequential Importance Sampling}
\begin{figure}[htbp]
\centering
	\includegraphics[trim=0cm 0cm 0cm 0cm,scale=0.25]{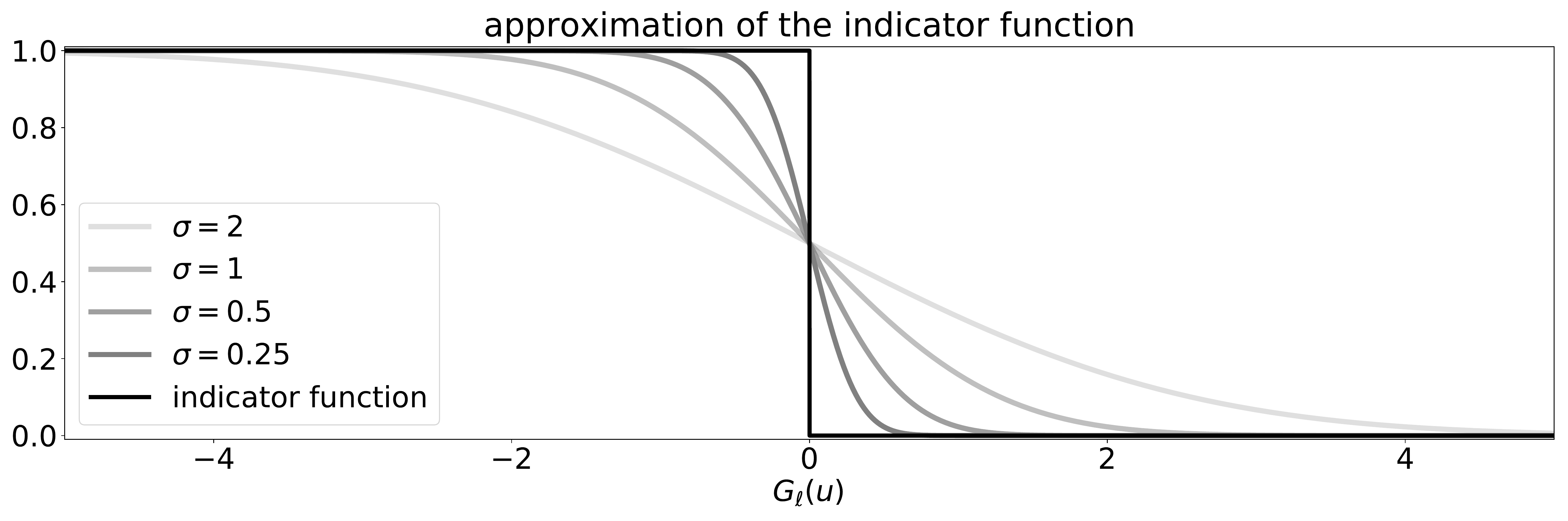}
		\caption[Indicator approximation]{Approximation of the indicator function $I(G_{\ell}(u)\le 0)$ by the cdf of the standard normal distribution $\Phi\left(-{G_{\ell}(u)}/{\sigma}\right)$.}
		\label{indicator approximation}
\end{figure}

According to \cite{Papaioannou16}, the sequence of IS densities is determined from a smooth approximation of the indicator function. The \emph{cumulative distribution function} (cdf) of the standard normal distribution is one possibility to approximate the indicator function. For $G_{\ell}(u)\neq 0$ we achieve pointwise convergence 
\begin{align*}
I(G_{\ell}(u)\le 0) = \lim_{\sigma \downarrow 0} \Phi\left(-\frac{G_{\ell}(u)}{\sigma}\right),
\end{align*}
while for $G_{\ell}(u)=0$ and $\forall \sigma>0$ it holds that $\Phi\left(-{G_{\ell}(u)}/{\sigma}\right)= 1/2 \neq I(G_{\ell}(u)\le 0)$, as visualized in Figure \ref{indicator approximation}. Further approximation functions are examined in \cite{Lacaze15} with an additional sensitivity analysis.
\\With the preceding consideration, the sequence of IS densities $\{p_{j,\ell}: j=0,\dots,N_T\}$ is defined as
\begin{align*}
p_{j,\ell}(u) &:= \frac{1}{P_{j,\ell}}\Phi\left(-\frac{G_{\ell}(u)}{\sigma_j}\right) \varphi_n(u) = \frac{1}{P_{j,\ell}}\eta_{j,\ell}(u), \text{ for } j=1,\dots, N_T,
\\ p_0(u) &:= \varphi_n(u),
\end{align*}
where $\infty > \sigma_1 > \cdots >\sigma_{N_T} > 0$ represent a strictly decreasing sequence of temperatures or bandwidths and $P_{j,\ell}$ is a normalizing constant such that $p_{j,\ell}$ is a well-defined density function. The denomination `temperatures' and their use is motivated by the temperature of the Boltzmann distribution \cite[Chapter VIII]{Gibbs02}. The number $N_T$ of tempering steps is a priori unknown and specifies the number of tempering steps to approximate the optimal IS density sufficiently accurately. Applying the IS approach, $P_{j,\ell}$ is determined by sampling from the density $p_{j-1,\ell}$ 
\begin{align}
P_{j,\ell} = \int_{\mathbb{R}^n} \eta_{j,\ell}(u) \mathrm{d}u = P_{j-1,\ell} \int_{\mathbb{R}^n} w_{j,\ell}(u) p_{j-1,\ell}(u) \mathrm{d}u = P_{j-1,\ell} \mathbb{E}_{p_{j-1,\ell}}[w_{j,\ell}(u)], \label{eq P_j}
\end{align}
where $w_{j,\ell}(u) := {\eta_{j,\ell}(u)}/{\eta_{j-1,\ell}(u)}$. Hence, the fraction of consecutive normalizing constants $S_{j,\ell}= {P_{j,\ell}}/{P_{j-1,\ell}}$ is estimated by
\begin{align}
\hat{S}_{j,\ell} := \hat{\mathrm{E}}_{p_{j-1,\ell}} [w_{j,\ell}(u)] = \frac{1}{N}\sum_{k=1}^{N} w_{j,\ell}(u_k),\label{S tempering}
\end{align}
where the samples $\{u_k\}_{k=1}^{N}$ are distributed according to $p_{j-1,\ell}$. Using the definition of $\eta_{j,\ell}$ and $\eta_{j-1,\ell}$, the weights $w_{j,\ell}(u_k)$ for $k=1,\dots,N$ are given by
\begin{align}
w_{j,\ell}(u_k) &= \frac{\Phi\left(-{G_{\ell}(u_k)}/{\sigma_j}\right)}{\Phi\left(-{G_{\ell}(u_k)}/{\sigma_{j-1}}\right)}, \text{ for } j>1\label{weights}
\\w_{1,\ell} (u_k) &= \Phi(-G_{\ell} (u_k) / \sigma_1).\notag
\end{align}
To obtain an accurate estimator $\hat{S}_{j,\ell}$, the parameters $\sigma_j$ are adaptively determined such that consecutive densities differ only slightly. This goal is achieved by requiring that the coefficient of variation of the weights $w_{j,\ell}$ is close to the target value $\delta_{\mathrm{target}}$, which is specified by the user. This leads to the following minimization problem
\begin{align}
\sigma_j = \underset{\sigma \in (0,\sigma_{j-1})}{\mathrm{argmin}} \big\Vert \delta_{w_{j,\ell}} - \delta_{\mathrm{target}}\big\Vert_2^2,\label{min pbl}
\end{align}
where $\delta_{w_{j,\ell}}$ is the \emph{coefficient of variation} of the weights (\ref{weights}). This adaptive procedure is similar to the adaptive tempering in \cite{Beskos13,Latz18} and is equivalent to requiring that the \emph{effective sample size} takes a target value \cite{Latz18}. Note that the solution of the minimization problem in (\ref{min pbl}) does not require further evaluations of the LSF. Hence, its costs are negligible compared to the overall computational costs. The tempering iteration is finished if the coefficient of variation $\delta_{w_{\mathrm{opt},\ell}}$ of the weights with respect to the optimal IS density
\begin{align}
w_{\mathrm{opt},\ell}(u_k) := I(G_{\ell}(u_k)\le 0)\frac{\varphi_n(u_k)}{\eta_{j,\ell}(u_k)}\label{weights optimal density}
\end{align}
is smaller than $\delta_{\mathrm{target}}$ and, hence, the optimal IS density is approximated sufficiently well. According to \cite{Papaioannou16}, the SIS estimator of the probability of failure is defined as follows
\begin{align}
\hat{P}_{f,\ell}^{\mathrm{SIS}} = \left(\prod_{j=1}^{N_T} \hat{S}_{j,\ell} \right)\frac{1}{N}\sum_{k=1}^{N} w_{\mathrm{opt},\ell}(u_k),\label{SIS estimator}
\end{align} 
where the weights $w_{\mathrm{opt},\ell}$ are defined in (\ref{weights optimal density}) with $j=N_T$. The sum over the weights $w_{\mathrm{opt},L}(u_k)$ in (\ref{SIS estimator}) represents the last tempering step from the IS density $p_{N_T,\ell}$ to the optimal IS density $p_{\mathrm{opt},\ell}$ given in (\ref{optimal IS density}). It corresponds to the estimator of the ratio $P_{N_T,\ell}/P_{f,\ell}$ since $P_{f,\ell}$ is the normalizing constant of the optimal IS density.
\\During the iteration, MCMC sampling is applied to transfer samples distributed according to $p_{j-1,\ell}$ to samples distributed according to $p_{j,\ell}$ for $j=1,\dots,N_T$. Section \ref{Section MCMC} explains MCMC sampling in more detail. Algorithm \ref{Tempering alg} summarizes the procedure of one tempering step for sampling from $p_{j,\ell}$ and estimating $S_{j,\ell}$ starting from samples from $p_{j-1,\ell}$. 

\begin{remark}
We remark that nestedness, which is a prerequisite for SuS, is not an issue for SIS. This is because the intermediate sampling densities are smooth approximations of the optimal IS density and they all have supports in the whole outcome space. The proximity of two consecutive densities is ensured by (\ref{min pbl}). This property of SIS motivates the development of MLSIS in the following section. We note that MLSuS does not satisfy nestedness, which leads to the denominators in the estimator (\ref{MLSuS estimator}). 
\end{remark}

\begin{algorithm}
\caption{Tempering algorithm ($N$ samples from $p_{j-1,\ell}$, $\sigma_{j-1}$, $\delta_{\mathrm{target}}$, $G_{\ell}$)}
\label{Tempering alg}
\begin{algorithmic}[1]
\STATE determine $\sigma_j$ from the optimization problem (\ref{min pbl})
\STATE evaluate the weights $w_{j,\ell}$ as in (\ref{weights}) for the current set of samples
\STATE evaluate the estimator $\hat{S}_{j,\ell}$ as in (\ref{S tempering})
\STATE re-sample the samples of $p_{j-1,\ell}$ based on their weights $w_{j,\ell}$
\STATE move the samples with MCMC to generate $N$ samples from the density $p_{j,\ell}$
\RETURN $N$ samples from $p_{j,\ell}$, $\sigma_j$, $\hat{S}_{j,\ell}$
\end{algorithmic}
\end{algorithm}

\section{Multilevel Sequential Importance Sampling}\label{Section MLSIS}
SIS and SuS have the drawback that all PDE solves are performed with the same discretization accuracy. This can lead to huge computational costs if the discretization level is high or the number of required tempering steps is large. Simply decreasing the level $\ell$ can lead to a bias in the estimated probability of failure, since the accuracy of the LSF decreases if the discretization level decreases. Therefore, the work in \cite{Latz18} develops the $\mathrm{MLS}^2\mathrm{MC}$ method, where computations are performed on a sequence of increasing discretization levels while achieving an improvement in terms of computational costs. Originally, this method has been developed for \emph{Bayesian inverse problems} \cite{Dashti2017}. In this section, we show how we can reformulate the $\mathrm{MLS}^2\mathrm{MC}$ method as an MLSIS estimator for the probability of failure.

\subsection{Bridging}
Consider the sequence of discretization levels $\ell\in\{1,\dots,L\}$ where $\ell=1$ represents the smallest and $\ell=L\in\mathbb{N}$ highest discretization level, i.e., finest element mesh size. Throughout this paper, it is assumed that the computational costs of evaluating $G_{\ell}$ are given by
\begin{align}
\mathrm{Cost}_{\ell} = \mathcal{O}(2^{-d(L-\ell)}),\label{Cost levels}
\end{align}
where $d\in\mathbb{N}$ is the dimension of the computational domain. In order to use a hierarchy of discretization levels, \emph{bridging} is applied to transfer samples following a distribution on a coarse grid to samples following a distribution on the next finer grid. The level update is defined as proposed in \cite{Koutsourelakis09}. The density $p_{j,\ell}$ of the coarse grid is transformed to the density $p_{j,\ell+1}$ of the next finer grid by the sequence
\begin{align}
p_{j,\ell}^{t}(u) := \frac{1}{P_{j,\ell}^{t}}\Phi\left(-\frac{G_{\ell+1}(u)}{\sigma_j}\right)^{\beta_{t}}\Phi\left(-\frac{G_{\ell}(u)}{\sigma_j}\right)^{1-\beta_{t}}\varphi_n(u),\label{eq bridging}
\end{align}
for $t=0,\dots,N_{B_{\ell}}$, where $0 = \beta_{0} < \beta_{1}< \cdots < \beta_{N_{B_{\ell}}}=1$, i.e., $p_{j,\ell}^0 = p_{j,\ell}$ and $p_{j,\ell}^{N_{B_{\ell}}}=p_{j,\ell+1}$. The number $N_{B_{\ell}}\in\mathbb{N}$ of intermediate bridging densities is a priori unknown. As in equation (\ref{eq P_j}), the quantity $P_{j,\ell}^t$ in (\ref{eq bridging}) can be calculated using samples distributed according to $p_{j,\ell}^{t-1}$. Similarly, the fraction of consecutive normalizing constants $S_{j,\ell}^t = P_{j,\ell}^t/P_{j-l}^{t+1}$ is estimated by
\begin{align}
\hat{S}_{j,\ell}^t:=\hat{\mathrm{E}}_{p_{j,\ell}^{t-1}} [w_{j,\ell}^t(u)] = \frac{1}{N}\sum_{k=1}^{N} w_{j,\ell}^t(u_k),\label{S bridging}
\end{align}
where the samples $\{u_k\}_{k=1}^N$ are distributed according to $p_{j,\ell}^{t-1}$ and the weights are given by
\begin{align}
w_{j,\ell}^t(u_k) := \frac{\Phi\left(-{G_{\ell+1}(u_k)}/{\sigma_j}\right)^{\beta_{t}}\Phi\left(-{G_{\ell}(u_k)}/{\sigma_j}\right)^{1-\beta_{t}}}{\Phi\left(-{G_{\ell+1}(u_k)}/{\sigma_j}\right)^{\beta_{t-1}}\Phi\left(-{G_{\ell}(u_k)}/{\sigma_j}\right)^{1-\beta_{t-1}}},\label{bridging weights}
\end{align}
for $k=1,\dots,N$. The bridging temperatures $\beta_t$ are adaptively determined by solving the minimization problem
\begin{align}
\beta_t = \underset{{\beta \in (\beta_{t-1}, 1]}}{\mathrm{argmin}} \big\Vert \delta_{w_{j,\ell}^t} - \delta_{\mathrm{target}} \Vert_2^2,\label{min bridging}
\end{align}
where $\delta_{w_{j,\ell}^t}$ is the coefficient of variation of the weights. As in \cite{Latz18}, we set the target coefficient of variation within the bridging steps to the same value as in the tempering steps. Within one level update, the bridging sequence is finished if $\beta_t = 1$ holds. Note that each level update requires a sequence of bridging densities and tempering is not performed during level updates. As in the tempering steps, MCMC sampling is applied to transfer samples between two consecutive bridging densities. By combining all estimators $\hat{S}$ of the tempering and bridging sequences given in (\ref{S tempering}) and (\ref{S bridging}), respectively, the MLSIS estimator for the probability of failure is given as
\begin{align}
\hat{P}_{f}^{\mathrm{MLSIS}} = \left(\prod_{j=1}^{N_T} \prod_{\ell=1}^{L} \prod_{t=1}^{N_{B_{\ell}}}\hat{S}_{j,\ell}^t \right)\frac{1}{N}\sum_{k=1}^N w_{\mathrm{opt},L}(u_k),\label{MLSIS estimator}
\end{align} 
where the weights $w_{\mathrm{opt},L}$ are defined in (\ref{weights optimal density}) with $j=N_T$ and represent the last tempering step from the IS density $p_{N_T,L}$ to the optimal IS density $p_{\mathrm{opt},L}$ given in (\ref{optimal IS density}). Algorithm \ref{Briging alg} summarizes the procedure of one level update.

\begin{algorithm}
\caption{Bridging algorithm ($N$ samples from $p_{j,\ell}$, $\sigma_{j}$, $\delta_{\mathrm{target}}$, $G_{\ell}$, $G_{\ell+1}$)}\label{Briging alg}
\begin{algorithmic}[1]
\STATE $t\leftarrow 0$
\STATE $\beta_t \leftarrow 0$
\WHILE{$\beta_t <1$}
	\STATE $t\leftarrow t+1$
	\STATE determine $\beta_{t}$ from the optimization problem (\ref{min bridging})
	\STATE evaluate the weights $w_{j,\ell}^t$ as in (\ref{bridging weights}) for the current set of samples
	\STATE evaluate the estimator $\hat{S}_{j,\ell}^t$ as in (\ref{S bridging})
	\STATE re-sample the samples of $p_{j,\ell}^{t-1}$ based on their weights $w_{j,\ell}^t$ 
	\STATE move the samples with MCMC to generate $N$ samples from the density $p_{j,\ell}^t$
\ENDWHILE
\RETURN $N$ samples from $p_{j,\ell+1}$, $\hat{S}_{j,\ell}^t$
\end{algorithmic}
\end{algorithm}

\subsection{Update scheme}\label{Update scheme}
The crucial part of the MLSIS method is to combine the adaptive tempering and bridging sequences and to provide a heuristic idea when to perform bridging or tempering. Initially, the samples $\{u_k\}_{k=1}^N$ are distributed according to the $n$-variate standard normal distribution $\varphi_n$, i.e., $\sigma_0 = \infty$. The LSF is evaluated on the smallest discretization level $\ell=1$. Tempering is always performed in the first step in order to determine $\sigma_1$ to approximate the indicator function. The tempering finishes if the coefficient of variation $\delta_{w_{\mathrm{opt},\ell}}$ of the weights with respect to the optimal IS density (\ref{weights optimal density}) is smaller than $\delta_{\mathrm{target}}$. The bridging finishes if the highest discretization level $\ell = L$ is reached. The combination of tempering and bridging determines costs and accuracy of the method. The authors in \cite{Latz18} analyse the efficiency of different decision schemes which leads to the following approach. The scheme should perform as many tempering steps as possible on small discretization levels while level updates are performed if the discrepancy between evaluations of two consecutive levels is too large. To measure this occurrence, a small subset of samples $\{u_{j_k}\}_{k=1}^{N_s}$ with $N_s<N$ is randomly selected without replacement. A level update is performed for this subset through one bridging step and the resulting coefficient of variation $\delta_{w^{N_s}}$ of the weights
\begin{align*}
w_{j,\ell}^{N_s}(u_{j_k}) = \frac{\Phi\left(-{G_{\ell+1}(u_{j_k})}/{\sigma_j}\right)}{\Phi\left(-{G_{\ell}(u_{j_k})}/{\sigma_j}\right)}, \text{ for } k=1,\dots,N_s,
\end{align*}
is estimated. Depending on the estimated value $\delta_{w^{N_s}}$, two cases occur:
\begin{itemize}
\item[1)] either $\delta_{w^{N_s}}> \delta_{\mathrm{target}}$: Bridging is performed since the accuracy is small, i.e., the difference between levels is high
\item[2)] or $\delta_{w^{N_s}} \leq \delta_{\mathrm{target}}$: Tempering is performed since the accuracy is high, i.e., the difference between levels is small.
\end{itemize}
If case 1) occurs, the evaluations with respect to $G_{\ell+1}$ can be stored and reused in the bridging step and invested costs are not wasted. Whereas in case 2), these evaluations are no longer required and invested costs are wasted. Calculating $\delta_{w^{N_s}}$ for the sample subset is redundant if tempering has already finished. Then, bridging is always performed to reach the final discretization level. Moreover as proposed in \cite{Latz18}, tempering is performed after each level update, if the tempering has not already finished. In this case, calculating $\delta_{w^{N_s}}$ is redundant, too. Note that $\delta_{w_{\mathrm{opt},\ell}}$ has to be calculated after each tempering and level update, to decide if tempering is finished. Finally, the MLSIS method is finished if both tempering and bridging are finished. The procedure is described in Algorithm \ref{MLSIS alg}.

\begin{remark}
We note that, according to \cite{Latz18}, the finest discretization level $L$ can be chosen adaptively based on the coefficient of variation $\delta_{w^{N_s}}$ between two consecutive discretization levels. Bridging is finished if $\delta_{w^{N_s}}$ is smaller than a given bound which is much smaller than $\delta_{\mathrm{target}}$.
\end{remark}

\begin{algorithm}
\caption{MLSIS algorithm ($N$, $n$, $L$, $\delta_{\mathrm{target}}$, $N_s$, $G_{\ell}$)}\label{MLSIS alg}
\begin{algorithmic}[1]
\STATE Generate $N$ samples from the $n$-variate standard normal distribution $\varphi_n$
\STATE $\ell \leftarrow 1$
\STATE Perform Tempering
\WHILE{Tempering is not finished \textbf{or} Bridging is not finished}
	\IF{Tempering is finished}
		\STATE Perform Bridging
		\STATE $\ell \leftarrow \ell + 1$
	\ELSIF{Bridging is finished \textbf{or} last step was a Bridging step}
		\STATE Perform Tempering
	\ELSE
		\STATE Perform Bridging in one step with a random subset of $N_s$ samples
		\STATE Calculate $\delta_{w^{N_s}}$
		\IF{$\delta_{w^{N_s}}<\delta_{\mathrm{target}}$}
			\STATE Perform Tempering
		\ELSE
			\STATE Perform Bridging
			\STATE $\ell \leftarrow \ell + 1$
		\ENDIF
	\ENDIF
	\STATE Calculate $\delta_{w_{\mathrm{opt},\ell}}$
	\IF{$\delta_{w_{\mathrm{opt},\ell}}\leq\delta_{\mathrm{target}}$}
		\STATE Tempering is finished
	\ENDIF
	\IF{$\ell=L$}
		\STATE Bridging is finished
	\ENDIF
\ENDWHILE
\RETURN Probability of failure estimate
\end{algorithmic}
\end{algorithm}

\subsection{Level dependent dimension}\label{level dependent dim}
As mentioned in Section \ref{SuS and MLSuS}, the nestedness problem of MLSuS motivates \cite{Ullmann15} to study a level dependent parameter dimension. This approach can also be applied in MLSIS to reduce variances between level updates and, hence, increase the number of tempering updates on coarse grids. For this purpose, it is assumed that the LSF $G$ depends on a random field that is approximated by a truncated \emph{Karhunen-Lo\`{e}ve} (KL) expansion. This setting occurs in many relevant applications as well as in numerical experiments presented in Section \ref{chapter numerical experiments}. Since high order KL terms are highly oscillating, they can not be accurately discretized on coarse grids which leads to noisy evaluations and higher variances. By reducing the number of KL terms on coarse grids, the variance between consecutive LSF evaluations is reduced. Therefore, the coefficient of variation $\delta_{w^{N_s}}$ is smaller and case 2) in Section \ref{Update scheme} is more likely. Hence, more tempering steps are performed on small discretization levels, which decreases the computational costs for MLSIS.

\section{Markov Chain Monte Carlo}\label{Section MCMC}
The goal of SIS and MLSIS is to transform samples from the prior density $p_0 = \varphi_n$ into samples of the optimal IS density $p_{\mathrm{opt}}$. Thereby, a sequence of densities is defined which converges to the optimal one. MCMC is applied to transform samples into consecutive densities of the tempering and bridging steps. 
\\Consider the tempering step from $p_{j-1,\ell}$ to $p_{j,\ell}$. The samples $\{u_k\}_{k=1}^{N}$ are distributed as $p_{j-1,\ell}$ and have to be transformed into samples that are distributed as $p_{j,\ell}$. To define the number of seeds of the MCMC algorithm, we choose a parameter 
\begin{align}
c\in(0,1] \text{ such that } \frac{1}{c}\in\mathbb{N} \text{ and } c\cdot N \in \mathbb{N}.\label{seed paramter}
\end{align} 
Then, $N_c:=c\cdot N$ seeds are randomly selected with replacement from the set $\{u_k\}_{k=1}^{N}$ according to their weights $\{w_{j-1,\ell}(u_k)\}_{k=1}^{N}$ given in (\ref{weights}). The set of seeds is denoted by $\{u_{k_j}\}_{j=1}^{N_c}$. In this procedure, which corresponds to \emph{multinomial resampling}, samples with high weights are copied multiple times and samples with low weights are discarded. There are also other resammpling methods, such as \emph{stratified resampling} or \emph{systematic resampling}, which can be applied. A study on their convergence behaviour is given in \cite{Gerber2019}. The burn-in length is denoted by $N_b\in\mathbb{N}$. Starting with the seed $u_0\in \{u_{k_j}\}_{j=1}^{N_c}$, a Markov chain of length $N_b+1/c$ is simulated that has $p_{j,\ell}$ as its stationary distribution. The first $N_b$ states are rejected after the simulation. Algorithm \ref{MCMC alg} states the MCMC procedure that employs the \emph{Metropolis-Hastings} sampler \cite{Hastings70,Metropolis53}. During the algorithm, a \emph{proposal} $\bar{u}$ is generated according to the \emph{proposal density} $q$. Moreover, the acceptance function $\alpha : \mathbb{R}^n\times \mathbb{R}^n\rightarrow [0,\infty)$ is given by
\begin{align*}
\alpha_T(u_0,\bar{u}) := \frac{\Phi\left(-{G(\bar{u})}/{\sigma_j}\right)\varphi_n(\bar{u})q(u_0\mid \bar{u})}{\Phi\left(-{G(u_0)}/{\sigma_j}\right)\varphi_n(u_0)q(\bar{u}\mid u_0)},
\end{align*}
which is the ratio of the target density $p_{j,\ell}$ with respect to the current state of the chain $u_0$ and candidate $\bar{u}$. For a bridging step, the seeds are selected from samples distributed according to $p_{j,\ell}^t$ and the target density is $p_{j,\ell}^{t+1}$. The weights are given by $\{w_{j,\ell}^t(u_k)\}_{k=1}^N$, see (\ref{bridging weights}), and the acceptance function $\alpha$ must be replaced by 
\begin{align*}
\alpha_B(u_0,\bar{u}) = \frac{\Phi\left(-{G_{\ell+1}(\bar{u})}/{\sigma_j}\right)^{\beta_{t+1}}\Phi\left(-{G_{\ell}(\bar{u})}/{\sigma_j}\right)^{1-\beta_{t+1}}\varphi_n(\bar{u})q(u_0\mid \bar{u})}{\Phi\left(-{G_{\ell+1}(u_0)}/{\sigma_j}\right)^{\beta_{t+1}}\Phi\left(-{G_{\ell}(u_0)}/{\sigma_j}\right)^{1-\beta_{t+1}}\varphi_n(u_0)q(\bar{u}\mid u_0)}.
\end{align*}
\begin{algorithm}
\caption{MCMC algorithm ($u_0$, $q(\cdot\mid \cdot)$, $\alpha(\cdot, \cdot)$, $c$, $N_b$)}\label{MCMC alg}
\begin{algorithmic}[1]
\STATE $\mathrm{Chain} = \emptyset$
\WHILE{$i\le N_b + 1/c$}
	\STATE Generate a candidate $\bar{u}$ from the proposal density $q(\cdot\mid u_0)$
	\STATE Evaluate $\alpha(u_0,\bar{u})$
	\STATE Accept the candidate $\bar{u}$ with probability $\min\{1,\alpha(u_0,\bar{u})\}$
	\IF{$\bar{u}$ is accepted}
		\STATE $u_0 \leftarrow \bar{u}$
	\ENDIF
	\STATE $\mathrm{Chain} \leftarrow \mathrm{Chain}\cup u_0$
	\STATE $i\leftarrow i+1$		
\ENDWHILE
\STATE Discard the first $N_b$ elements of $\mathrm{Chain}$
\RETURN simulation of Markov chain
\end{algorithmic}
\end{algorithm} 
\begin{remark}
Since consecutive densities within SIS and MLSIS are constructed in a way that they are not too different and samples are weighted according to the target distribution, the burn-in length can be small or even negligible within SIS and MLSIS \cite{Papaioannou16}. Note that for SuS and MLSuS, the $N\cdot \hat{p}_0$ samples with the lowest LSF values are selected as seeds.
\end{remark}

\subsection{Adaptive conditional sampling}
The random walk Metropolis Hastings algorithm \cite{Hastings70,Metropolis53} is a classical MCMC algorithm. However, random walk samplers suffer from the curse of dimensionality, i.e., the acceptance rate is small in high dimensions, see \cite{Papaioannou15}. Since high dimensional parameter spaces are considered in the numerical experiments in Section \ref{chapter numerical experiments}, \emph{adaptive conditional sampling} (aCS) is proposed, where the chain correlation is adapted to ensure a high acceptance rate. aCS is a dependent MCMC algorithm, i.e., the proposal density depends on the current seed $u_0$. More formally, the proposal $q$ is defined as the conditional multivariate normal density with mean vector $\rho u_0$ and covariance matrix $\Sigma=(1-\rho^2) I_n$, with $I_n$ denoting the identity matrix. During the iterations, $\rho\in [0,1]$ is adaptively adjusted such that the acceptance rate is around $44\%$ \cite{roberts2001}. This value leads to an optimal value in terms of the minimum autocorrelation criterion. By the structure of the proposal, the acceptance functions read as 
\begin{align*}
\alpha_T(u_0,\bar{u}) &= \frac{\Phi\left(-{G(\bar{u})}/{\sigma_j}\right)}{\Phi\left(-{G(u_0)}/{\sigma_j}\right)},
\\\alpha_B(u_0,\bar{u}) &= \frac{\Phi\left(-{G_{\ell+1}(\bar{u})}/{\sigma_j}\right)^{\beta_{t+1}}\Phi\left(-{G_{\ell}(\bar{u})}/{\sigma_j}\right)^{1-\beta_{t+1}}}{\Phi\left(-{G_{\ell+1}(u_0)}/{\sigma_j}\right)^{\beta_{t+1}}\Phi\left(-{G_{\ell}(u_0)}/{\sigma_j}\right)^{1-\beta_{t+1}}}.
\end{align*}
A more detailed description of the algorithm with the adaptive adjustment of the correlation parameter is given in \cite{Papaioannou15}. The aCS algorithm can be viewed as an adaptive version of the \emph{preconditioned Crank-Nicolson} (pCN) sampler \cite{Cotter13} tailored to application within SIS. 

\subsection{Independent sampler with von Mises-Fisher Nakagami proposal distribution}
\begin{figure}[htbp]
	\centering
	\includegraphics[trim=0cm 0cm 0cm 0cm,scale=0.23]{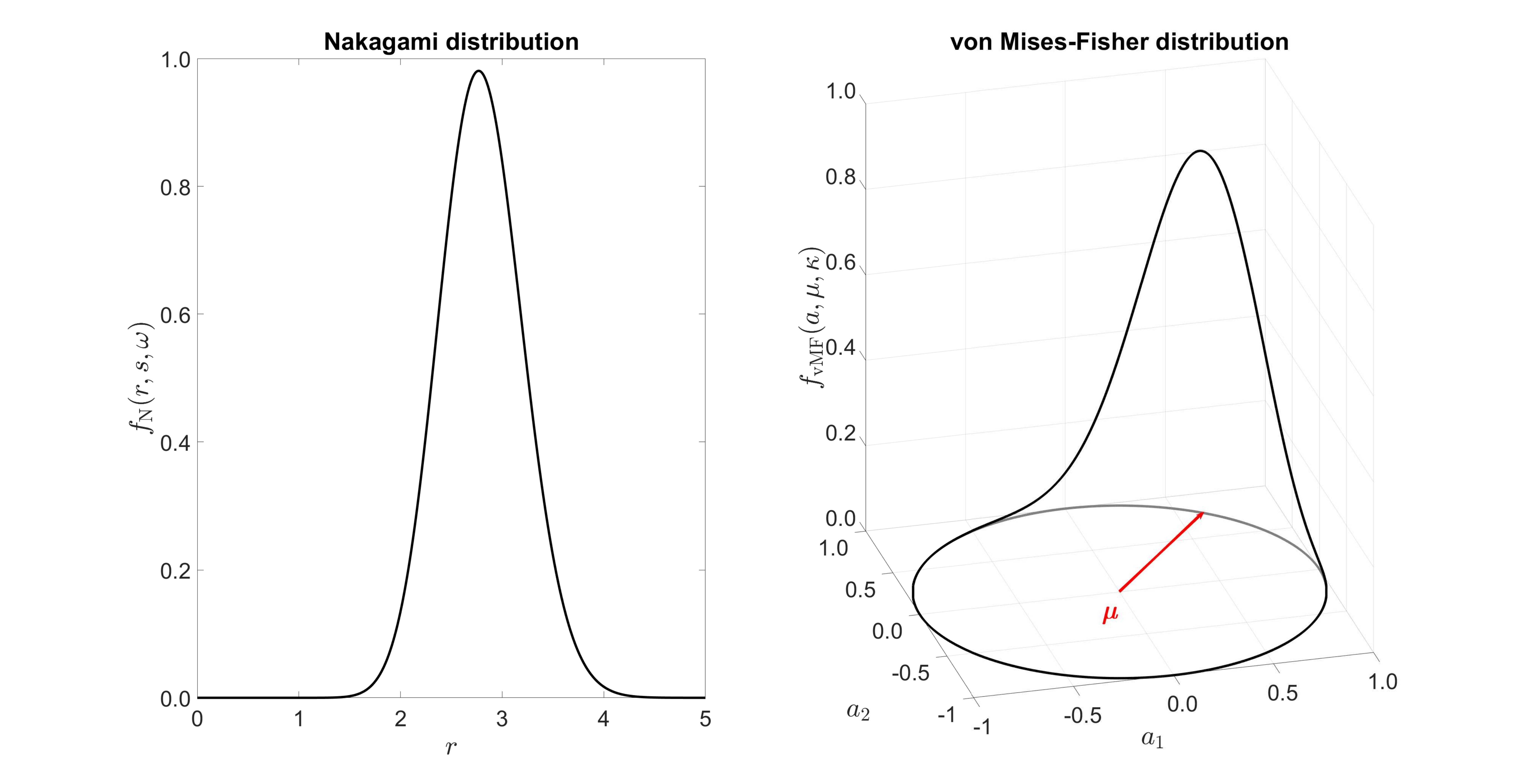}
		\caption[vMFN]{Illustration of the Nakagami distribution (left) and von Mises-Fisher distribution (right) in two dimensions. The parameters are defined as $\nu= (0.6, 0.75)^T$, $\kappa = 11$, $s=12$ and $\gamma=8$.}
		\label{Fig: vMFN}
\end{figure}

Since aCS is a dependent MCMC algorithm, the states of the chains are correlated which can lead to a higher variance of the estimated ratio of normalizing constants $\hat{S}_{j,\ell}$ and $\hat{S}_{j,\ell}^t$ given in (\ref{S tempering}) and (\ref{S bridging}), respectively. Hence, this leads to a higher variance of the estimated probability of failure (\ref{MLSIS estimator}). An independent MCMC algorithm overcomes this problem through using a proposal density that does not depend on the current state. The dependence on the current state enters in the acceptance probability. If the proposal density is chosen close to the target density, the acceptance probability will be close to one and the samples will be approximately independent. In the context of SIS and MLSIS, the available samples and corresponding weights of each previous density can be used to fit a distribution model to be used as proposal density in the MCMC step \cite{Chopin01,Papaioannou16}. For instance, \emph{Gaussian mixture} models can be used as a proposal density \cite{Papaioannou16}. A drawback of Gaussian densities in high dimensions is the concentration of measure around the hypersphere with norm equal to $\sqrt{n}$, see \cite{Katafygiotis08,Papaioannou19}. Therefore, only the direction of the samples is of importance. Furthermore, the Gaussian mixture model with $K$ densities has $K {n(n+3)}/{2}+(K-1)$ parameters, which have to be estimated. Both issues motivate the vMFN distribution. Therein, the direction is sampled from the \emph{von Mises-Fisher} (vMF) distribution \cite{Wang16} while the radius is sampled from the \emph{Nakagami} distribution, see \cite{Nakagami60}. The vMFN mixture model has only $K(n+3)+(K-1)$ parameters, which scales linearly in the dimension $n$. Note that for the Gaussian mixture, the number of parameters of the distribution model scales quadratically in the dimension $n$. To apply the vMFN distribution as the proposal density in Algorithm \ref{MCMC alg}, the parameters of the distribution model have to be fitted in advance.
\\It is assumed that all samples $u\in\mathbb{R}^n$ are given in their polar coordinate representation $u=r\cdot a$, where $r=\Vert u\Vert_2\in\mathbb{R}_+$ is the norm of $u$ and $a={u}/{\Vert u \Vert_2}\in\mathbb{R}^n$ its direction. For $u= r\cdot a$ the vMFN distribution with one mixture is defined as the product of the von Mises-Fisher and the Nakagami distribution, that is,
\begin{align*}
f_{\mathrm{vMFN}}(r,a\mid\nu,\kappa,s,\gamma) = f_{\mathrm{N}}(r\mid s,\gamma)\cdot f_{\mathrm{vMF}}(a\mid\nu,\kappa).
\end{align*}
The vMF distribution $f_{\mathrm{vMF}}$ defines the distribution of the direction on the $n$-dimen\-sional hypersphere $\mathbb{S}^{n-1}:=\{x\in\mathbb{R}^n: \Vert x \Vert_2 = 1\}$ and is given by
\begin{align*}
f_{\mathrm{vMF}}(a\mid\nu,\kappa) = \frac{\kappa^{n/2-1}}{(2\pi)^{n/2}\mathcal{I}_{n/2-1}(\kappa)}\exp(\kappa\nu^Ta),
\end{align*}
where $\nu\in\mathbb{S}^{n-1}$ is a mean direction and $\kappa\ge 0$ characterises the concentration around $\nu$. $\mathcal{I}_{n/2-1}$ denotes the \emph{modified Bessel function} of the first kind and order $n/2-1$ \cite[Chapter 9]{Abramowitz64}. Contrary, the Nakagami distribution $f_\mathrm{N}$ specifies the distribution of the radius and is defined by
\begin{align*}
f_{\mathrm{N}}(r\mid s,\gamma) := \frac{2s^s}{\Gamma(s)\gamma^s}r^{2s-1} \exp\left(-\frac{s}{\gamma}r^2\right),
\end{align*}
where $\Gamma(s)$ is the \emph{Gamma function}, $s\ge0.5$ is a shape parameter and $\gamma>0$ a spread parameter. Figure \ref{Fig: vMFN} shows an illustration of $f_{\mathrm{N}}$ and $f_{\mathrm{vMF}}$ for certain parameter values.
\begin{remark}
We have defined the vMFN distribution for the radius and direction $(r,a)$ on $[0,\infty)\times \mathbb{S}^{n-1}$. However, we actually approximate a distribution on $\mathbb{R}^n$, which defines the distribution of $u=r\cdot a\in\mathbb{R}^n$. By \cite[Theorem 1.101]{Klenke2014}, the distribution of $u\in\mathbb{R}^n$, which is the product distribution of $r$ and $a$, is given as 
\begin{align*}
f_U(u) &= \int_0^{\infty} f_{\mathrm{N}}(r\mid s,\gamma)f_{\mathrm{vMF}}\left(\frac{u}{r}\mid \nu,\kappa\right)\frac{1}{r^n}\mathrm{d}r
\\&\propto\int_0^{\infty}r^{2s-1-n}\exp\left(-\frac{s}{\gamma}r^2+\frac{\kappa\nu^T u}{r}\right)\mathrm{d}r.
\end{align*}
Since we can easily separate $u$ into $r$ and $a$, we usually work with $f_{\mathrm{vMFN}}$ rather than $f_U$.
\end{remark}
To define the vMFN distribution as a proposal density for Algorithm \ref{MCMC alg}, the parameters $\nu,\kappa,s$ and $\gamma$ have to be fitted using the current set of samples $\{u_k=r_k\cdot a_k\}_{k=1}^{N}$ and their weights $\{w_k\}_{k=1}^{N}$, which are given by (\ref{weights}) or (\ref{bridging weights}) for a tempering or bridging step, respectively. The parameters are determined by maximizing the weighted log-likelihood
\begin{align*}
\max_{\nu, \kappa, s, \gamma} \sum_{k=1}^N w_k \ln(f_{\mathrm{vMFN}}(r_k,a_k\mid \nu, \kappa, s, \gamma)).
\end{align*}
Differentiating this expression with respect to the parameters and setting the derivatives equal to zero yields the optimal parameters for the fitting \cite{Papaioannou19}. However, for the concentration $\kappa$ and shape parameter $s$ we use an approximation since the derivatives require the solutions of non-linear equations which arise from the Gamma function and modified Bessel function \cite{Bouhlel15,Wang16}. The fitted mean direction $\hat{\nu}$ and concentration $\hat{\kappa}$ are given by
\begin{align}
\hat{\nu} = \frac{\sum_{k=1}^N w_k\cdot a_k}{\Vert \sum_{k=1}^N w_k\cdot a_k\Vert_2},\hspace{0.5cm} \hat{\kappa} = \frac{\chi\cdot n - \chi^3}{1-\chi^2},\text{ where } \chi = \min\left\lbrace\frac{\Vert \sum_{k=1}^N w_k\cdot a_k\Vert_2}{\sum_{k=1}^N w_k}, 0.95\right\rbrace.\label{vMFN nu calc}
\end{align}
The upper bound of $0.95$ in (\ref{vMFN nu calc}) is chosen to ensure numerical stability of the algorithm. If $\chi$ converges to $1$, the vMFN distribution would converge to a point density \cite{Papaioannou19}. Moreover for the Nakagami distribution, the fitted spread $\hat{\gamma}$ and shape parameter $\hat{s}$ are given by
\begin{align*}
\hat{\gamma} = \frac{\sum_{k=1}^N w_k \cdot r_k^2}{\sum_{k=1}^N w_k}, \hspace{0.5cm} \hat{s} = \frac{\hat{\gamma}^2}{\nu_4-\hat{\gamma}^2}, \text{ where } \nu_4 = \frac{\sum_{k=1}^N w_k\cdot r_k^4}{\sum_{k=1}^N w_k}.
\end{align*}
To apply Algorithm \ref{MCMC alg} with respect to the vMFN distribution, the proposal $q(\cdot \mid u_0)$ is replaced by $f_{\mathrm{vMFN}}(\cdot,\cdot,\hat{\nu},\hat{\kappa},\hat{s},\hat{\gamma})$ with the fitted parameters.
\begin{remark}
If a mixture of vMFN distributions is considered with $K>1$ individual vMFN densities, the vMFN mixture distribution reads as 
\begin{align*}
f_{\mathrm{vMFNM}}(r,a\mid\boldsymbol{\nu},\boldsymbol{\kappa},\boldsymbol{s},\boldsymbol{\gamma}) = \sum_{j=1}^K \pi_j f_{\mathrm{vMFN}}(r,a\mid\nu_j,\kappa_j,s_j,\gamma_j),
\end{align*}
where the weights $\pi_{j}$ represent the probability of each mode and $\sum_{j=1}^K \pi_j=1$. In this case, the assignments of the samples to the modes is unknown and this assignment has to be estimated in addition. Therefore, the required parameters cannot be estimated in one iteration. For instance, the \emph{Expectation-Maximization} algorithm \cite{McLachlan05} can be applied to estimate the parameters iteratively. The resulting formulas are given in \cite{Papaioannou19}. The usage of mixtures is motivated by multi-modal failure domains. In the numerical experiments in Section \ref{chapter numerical experiments} only $K=1$ is considered.
\end{remark}

\subsection{MCMC for a level dependent dimension}
In the case that MLSIS or MLSuS are applied with a level dependent parameter dimension, the procedure of a level update has to be adjusted. Consider the level update from level $\ell$ to $\ell + 1$ and assume that the corresponding LSFs are defined as $G_\ell : \mathbb{R}^{n_{\ell}}\rightarrow \mathbb{R}$ and $G_{\ell+1} : \mathbb{R}^{n_{\ell+1}} \rightarrow \mathbb{R}$, respectively, where $n_{\ell} < n_{\ell+1}$. Before the first MCMC step of the level update is carried out, the weights $w_{j,\ell}^1(u_k)$ for $k=1,\dots,N$ (see (\ref{bridging weights})) have to be evaluated. However, these evaluations require the evaluation of $G_{\ell+1}$ with respect to the current samples $\{u_{k}\}_{k=1}^{N}$ which are defined on $\mathbb{R}^{n_{\ell}}$. In the beginning of MLSIS or MLSuS the samples $u_k$ are initialized from the standard normal density $\varphi_{n_{1}}$. Therefore, it is natural to sample the missing dimensions $\Delta n_{\ell+1} = n_{\ell+1}-n_{\ell}$ from the standard normal distribution $\varphi_{\Delta n_{\ell+1}}$. Hence, for each $k=1,\dots, N$ we sample $\Delta n_{\ell+1}$ independent standard normal random variables $\psi_k\in\mathbb{R}^{\Delta n_{\ell+1}}$ and stack $u_k$ and $\psi_k$ together, i.e., $\tilde{u}_k = [u_k, \psi_k]\in\mathbb{R}^{n_{\ell+1}}$. In order to evaluate the weights $w_{j,\ell}^1(u_k)$, the LSF $G_{\ell+1}$ is evaluated for $\tilde{u}_k$ and $G_{\ell}$ for $u_k$. The seeds for the MCMC step are chosen based on these weights. Subsequently, Algorithm \ref{MCMC alg} is performed. Within the MCMC algorithm, a proposal $\bar{u}\in\mathbb{R}^{n_{\ell+1}}$ is sampled from $q(\cdot\mid u_0)$ which is suitable for the evaluations of $G_{\ell+1}$. For the LSF $G_{\ell}$ the first $n_{\ell}$ entries of $\bar{u}$ are taken as input.

\section{Numerical experiments}\label{chapter numerical experiments}
In the following examples, all probability of failure estimates are obtained with respect to the same, finest discretization level, i.e., the multilevel methods iterate until this level is reached and the single-level methods estimator are based on this level. Therefore, the obtained errors involve only sampling errors while discretization errors are not included. 

\subsection{1D diffusion equation}
\begin{figure}[htbp]
\centering
	\includegraphics[trim=0cm 0cm 0cm 0cm,scale=0.23]{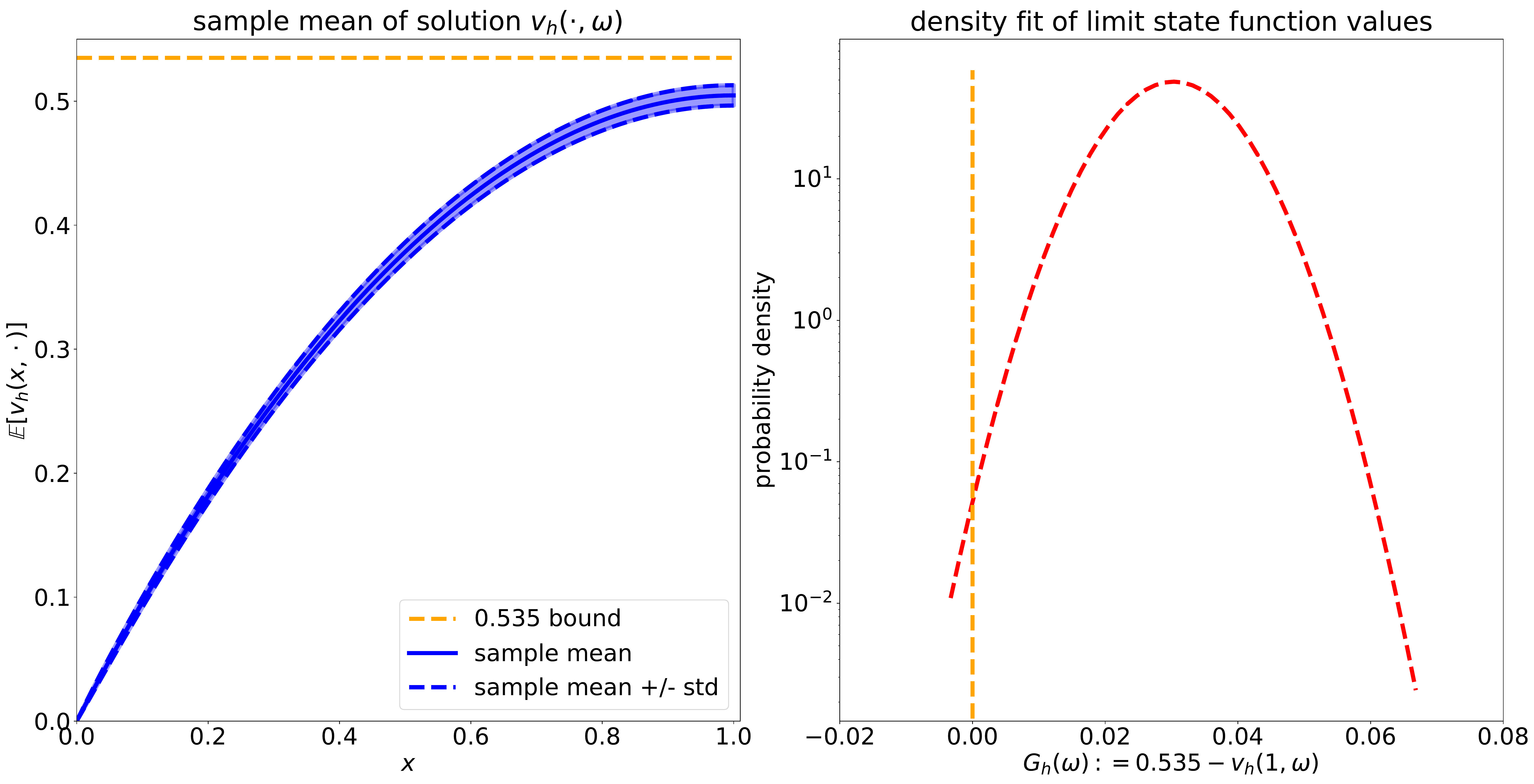}
		\caption[1D histogram]{Mean over $10^5$ realizations of solutions $v_h(\cdot,\omega)$ plus/minus the standard deviation (left). Fit of the respective LSF values $G_h(\omega)$ for $h={1}/{512}$ (right). Note that the probabilities are shown on a log-scale.}
		\label{histogram}
\end{figure}

We begin with Example 2 in \cite{Ullmann15} which considers the diffusion equation in the one-dimensional domain $D=[0,1]$. In particular, the corresponding stochastic differential equation is given by
\begin{align}
-\frac{\partial}{\partial x}\left(a(x, \omega)\frac{\partial}{\partial x} v(x,\omega)\right) &= 1 \text{ for } 0\le x \le 1,\label{1D diffusion}
\\ \text{such that} \hspace{0.2cm} v(0,\omega) &= 0 \text{ and } v'(1,\omega) = 0,\notag
\end{align}
for almost every (a.e.) $\omega\in\Omega$. Failure is defined as the event that the solution $v$ is larger than $0.535$ at $x=1$, i.e., $G(\omega):= 0.535 - v(1,\omega)\le 0$. The solution $v$ is approximated by a piecewise linear, continuous FEM approximation $v_h$ on a uniform grid with mesh size $h>0$. Hence, the approximated LSF is given by $G_{\ell}(\omega) = 0.535-v_{h_\ell}(1,\omega)$, where $\ell\in\mathbb{N}$ defines the discretization level. By crude Monte Carlo sampling (\ref{MC esimator}) with $N=10^7$ samples on a grid with mesh size $h={1}/{512}$, the probability of failure is estimated to be $P_f = 1.524\cdot 10^{-4}$. In the following, this value is referred to as the reference solution. Figure \ref{histogram} shows the mean of $10^5$ realizations of solutions $v_h(\cdot,\omega)$ plus/minus the standard deviation for $h={1}/{512}$. Additionally, the respective histogram of their LSF values is presented. We see that very few realizations are larger than $0.535$ at $x=1$.
\\The coefficient function $a(x,\omega)=\exp(Z(x,\omega))$ in (\ref{1D diffusion}) is a log-normal random field with constant mean function $\mathrm{E}[a(x,\cdot)] = 1$ and standard deviation $\mathrm{Std}[a(x,\cdot)] = 0.1$. That is, $Z$ is a Gaussian random field with constant mean function $\mu_Z = \log(\mathrm{E}[a(x,\cdot)]) - {\zeta_Z^2}/{2}$ and variance $\zeta_Z^2 = \log\left(({\mathrm{Std}[a(x,\cdot)]^2 + \mathrm{E}[a(x,\cdot)]^2})/{\mathrm{E}[a(x,\cdot)]^2}\right).$ Moreover, $Z$ has an exponential type covariance function which is given by $c(x,y) = \zeta_Z^2\exp\left(-{\vert x -y\vert}/{\lambda}\right)$, where $\lambda = 0.01$ denotes the correlation length. The infinite-dimensional log-normal random field $a$ is discretized by the truncated KL expansion of $Z$
\begin{align*}
Z(x,\omega) = \mu_Z + \zeta_Z \sum_{m=1}^M \sqrt{\nu_m}\theta_m(x) U_m(\omega),
\end{align*}
where $(\nu_m, \theta_m)$ are the KL eigenpairs and $\{U_m\}_{m=1}^M$ are independent standard normal Gaussian random variables. The eigenpairs can be analytically calculated as explained in \cite[p. 26ff]{Ghanem91}.
\\The probability of failure is estimated by SIS, MLSIS, SuS and MLSuS. For all methods, the estimation is performed for $N=250, 500, 1000, 2000$ samples and $N_s = 0.1\cdot N$ samples are considered for the small sample subset to decide if either bridging or tempering is performed in the update scheme of Section \ref{Update scheme}. For each parameter setting, the estimation is repeated $100$ times. For the multilevel methods, the sequence of mesh sizes is $h_\ell = 2^{-\ell-1}$ for $\ell=1,..,8$, i.e., the coarsest mesh size is $h_1={1}/{4}$ and the finest $h_8={1}/{512}$. If a level dependent dimension is considered, the parameter dimensions of the KL expansions are $n_1=10$, $n_2=20$, $n_3=40$, $n_4=80$ and $n_5=n_6=n_7=n_8=150$ as proposed in \cite{Ullmann15}. For a fixed parameter dimension, the dimension is $n=150$ for all discretization levels. This captures $87\%$ of the variability of $\log(a)$ \cite{Ullmann15}. SIS and MLSIS are performed for target coefficient of variations $\delta_{\mathrm{target}}=0.25$ and $\delta_{\mathrm{target}}=0.50$, which is considered in (\ref{min pbl}), (\ref{min bridging}). aCS and the independent sampler with the vMFN distribution and one mixture are considered as the MCMC methods without a burn-in. The parameter $c$ to define the number of seeds of the MCMC algorithm in (\ref{seed paramter}) is $c=0.1$ or $c=1$. For SuS and MLSuS, aCS is considered as the MCMC method without a burn-in and the parameter $\hat{p}_0$ in (\ref{SuS parameter}) is $\hat{p}_0 =0.1$ or $\hat{p}_0 =0.25$.

\subsubsection{Results}
\begin{figure}[h!]
\centering
	\includegraphics[trim=0cm 0cm 0cm 0cm,scale=0.23]{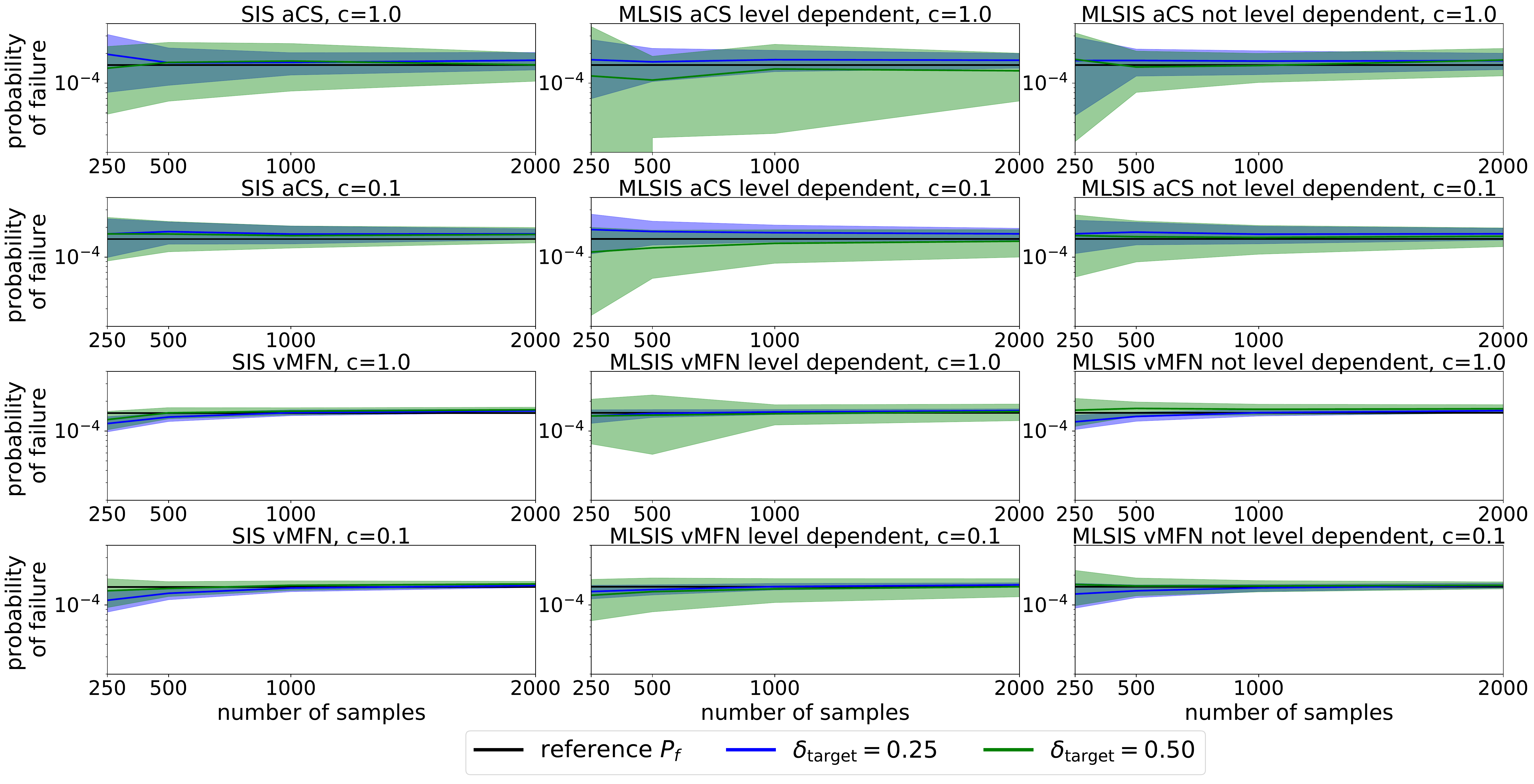}
		\caption[a]{Estimated probability of failure by SIS and MLSIS averaged over $100$ runs for $250, 500, 1000$ and $2000$ samples and $\delta_{\mathrm{target}}\in\{0.25, 0.50\}$. The coloured areas show the standard deviation of the estimates. The black lines show the reference estimate by Monte Carlo sampling. 1st row: aCS, $c=1.0$; 2nd row: aCS, $c=0.1$; 3rd row: vMFN, $c=1.0$; 4th row: vMFN, $c=0.1$; 1st column: SIS; 2nd column: MLSIS with level dependent dimension; 3rd column: MLSIS without level dependent dimension.}
		\label{1D-Diffusion: probability of failure + standard deviation}
\end{figure}

\begin{figure}[h!]
\centering
	\includegraphics[trim=0cm 0cm 0cm 0cm,scale=0.23]{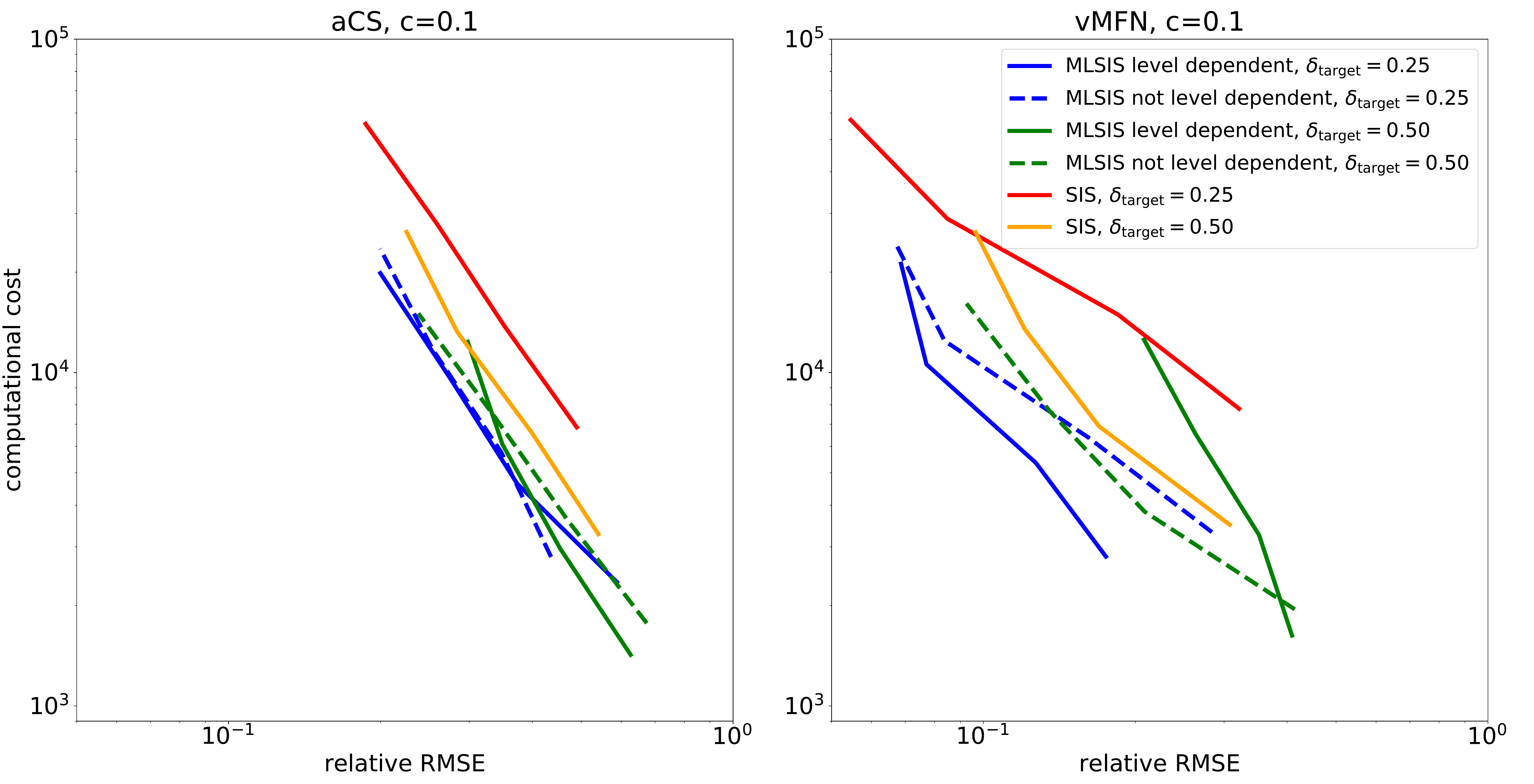}
		\caption[1D-Diffusion: RMSE and Costs]{Computational costs and relative RMSE of SIS and MLSIS averaged over $100$ runs for $250, 500, 1000$ and $2000$ samples and $\delta_{\mathrm{target}}\in\{0.25, 0.50\}$. 1st column: aCS, $c=0.1$; 2nd column: vMFN, $c=0.1$.}
		\label{1D-Diffusion: RMSE and Costs}
\end{figure}

\begin{figure}[h!]
	\centering
	\includegraphics[trim=0cm 0cm 0cm 0cm,scale=0.23]{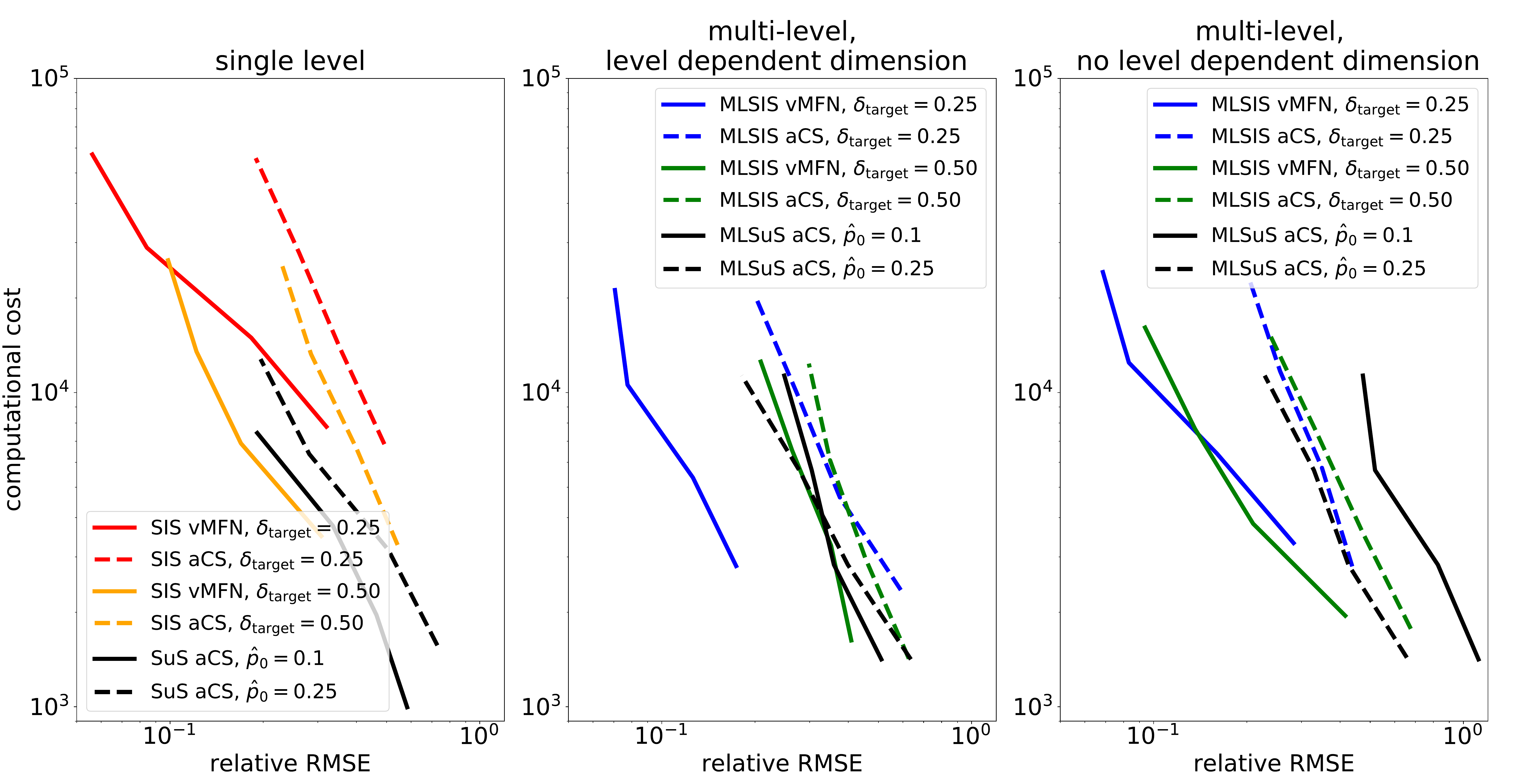}
		\caption[1D-Diffusion: Comparison RMSE and Costs]{Computational costs and relative RMSE of SIS, MLSIS, SuS and MLSuS averaged over $100$ runs for $250, 500, 1000$ and $2000$ samples. SIS and MLSIS are considered with aCS and vMFN, $c=0.1$ and $\delta_{\mathrm{target}}\in\{0.25, 0.50\}$. SuS and MLSuS are considered with aCS and $\hat{p}_0 =0.1$ or $\hat{p}_0 =0.25$. 1st column: single level; 2nd column: multi-level with level dependent dimension; 3rd column: multi-level without level dependent dimension.}
		\label{1D-Diffusion: Comparison RMSE and Costs}
\end{figure}

Figure \ref{1D-Diffusion: probability of failure + standard deviation} shows the estimated mean probability of failure by SIS and MLSIS plus/minus its standard deviation. The estimates of the means are in accordance with the reference solution for all settings. As expected, the bias and standard deviation decrease with an increasing number of samples. Furthermore, the standard deviation is smaller for a smaller target coefficient of variation. We observe that sampling from the vMFN distribution with independent MCMC yields a smaller bias and smaller standard deviation than applying aCS. Additionally, we observe that $c=0.1$ yields also a smaller standard deviation than $c=1$. Comparing the MLSIS results with the SIS results for $\delta_{\mathrm{target}}=0.50$, we see that SIS reaches a smaller standard deviation. For $\delta_{\mathrm{target}}=0.25$ the results are similar. For MLSIS and $\delta_{\mathrm{target}}=0.50$, a level dependent parameter dimension leads to a higher standard deviation than a fixed parameter dimension. However, for $\delta_{\mathrm{target}}=0.25$, the results are similar. We summarize that $\delta_{\mathrm{target}}=0.25$ yields for all settings a similar bias and standard deviation. Only the MCMC algorithm has a larger influence on the standard deviation in this setting. 
\\Figure \ref{1D-Diffusion: RMSE and Costs} shows the \emph{relative root mean square error} (RMSE) on the horizontal axis and the computational costs on the vertical axis of the SIS and MLSIS estimators. The relative RMSE is defined as
\begin{align*}
\mathrm{relRMSE} := \frac{\left(\mathbb{E}\Bigl\lbrack\left(\hat{P}_f-P_f\right)^2\Bigr\rbrack\right)^{\frac{1}{2}}}{P_f},
\end{align*}
where $\hat{P}_f$ denotes the estimated probability of failure. The costs are calculated based on the formula given in (\ref{Cost levels}) for $L=8$ and $d=1$. SIS and MLSIS yield a similar range of the relative RMSE but the computational costs are lower for MLSIS. Comparing the computational costs shown in Figure \ref{1D-Diffusion: RMSE and Costs}, we can save around $61\%$ of the computational costs if we apply MLSIS for the estimation. This shows the achievement of the main goal of the MLSIS algorithm, that is to save computational costs by employing a hierarchy of discretization levels. Furthermore, we observe that sampling from the vMFN distribution yields a lower relative RMSE than applying aCS. In case of $\delta_{\mathrm{target}}=0.25$, a level dependent dimension yields a smaller relative RMSE and lower computational costs than a fixed parameter dimension. This was expected since variances between level updates are smaller and, therefore, more tempering steps are performed on coarse levels. However, MLSIS with a level dependent dimension, $\delta_{\mathrm{target}}=0.50$ and sampling from the vMFN distribution yields a higher relative RMSE than applying a fixed parameter dimension for the same computational cost.
\\Figure \ref{1D-Diffusion: Comparison RMSE and Costs} shows the relative RMSE and computational costs of SIS, MLSIS, SuS and MLSuS. We observe that SuS yields the same relative RMSE as SIS with aCS. However, SuS requires less computational costs. If we consider SIS with vMFN, the relative RMSE is smaller compared to SuS but the computational costs are higher for SIS. For the multilevel methods with a level dependent dimension, we observe that MLSuS and MLSIS with aCS yield a similar relative RMSE but MLSIS requires more computational costs. However, the savings with MLSuS are smaller compared to the single-level estimators. MLSIS with vMFN and $\delta_{\mathrm{target}}=0.25$ yields a much smaller relative RMSE than all other estimators and computational costs can be saved compared to MLSuS. Theses results are similar to the multilevel methods without a level dependent dimension. In this case, we can observe that MLSuS with $\hat{p}_0=0.1$ yields a large relative RMSE which is due to the nestedness issue of MLSuS.

\subsection{2D flow cell}
\begin{figure}[h!]
\centering
	\includegraphics[trim=0cm 0cm 0cm 0cm,scale=0.3]{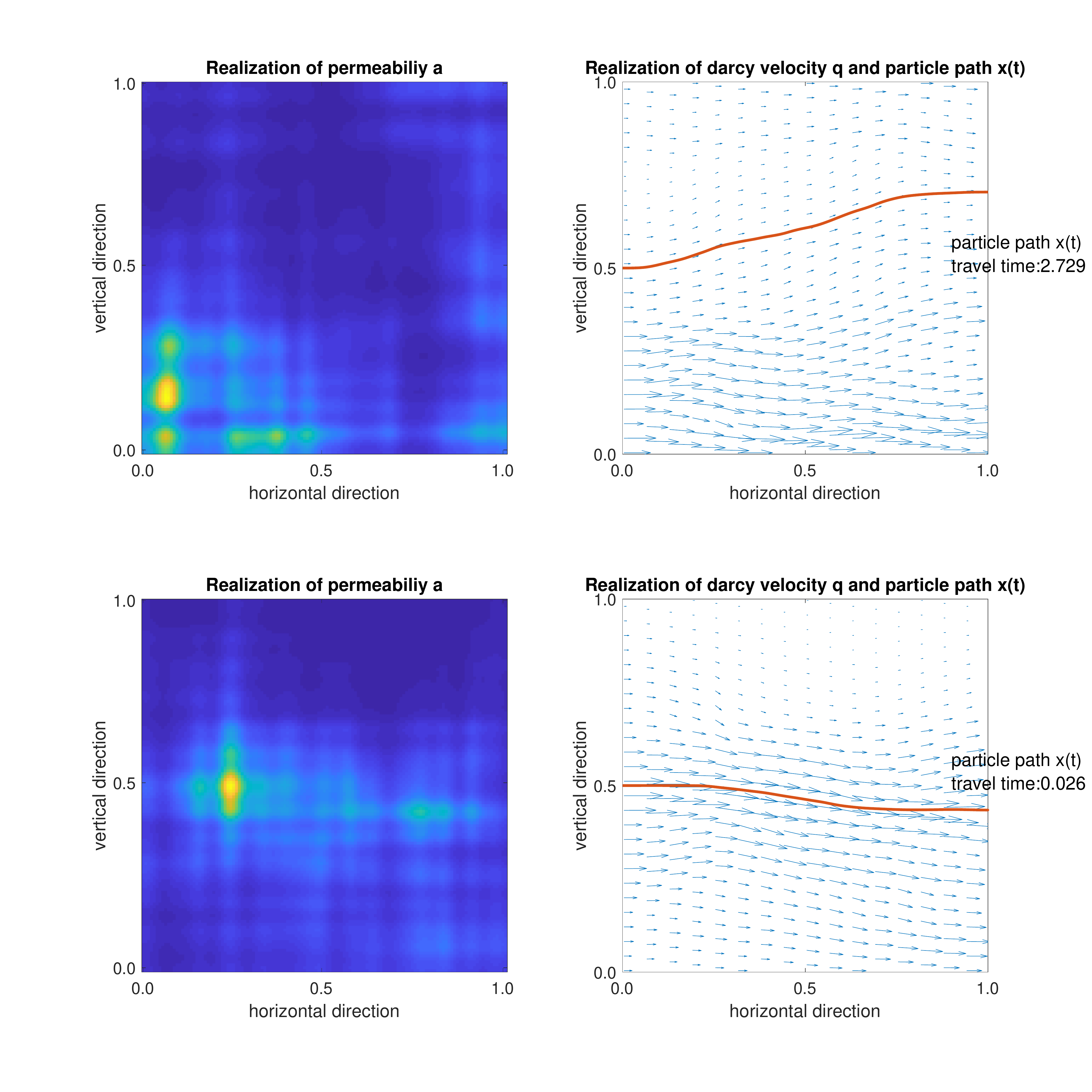}
		\caption[flow cell realization]{Realization of the permeability $a(\cdot,\omega)$ and the respective solution of the Darcy velocity $q_h(\cdot,\omega)$ and the particle path $x(t,\omega)$ for $h={1}/{128}$.}
		\label{flow_cell}
\end{figure}

We consider the two-dimensional application in \cite[Section 6.1]{Ullmann15}, which is a simplified setting of the rare event arising in planning a radioactive waste repository, see Section \ref{Sec: Introduction}. The probability of failure is based on the travel time of a particle within a two-dimensional flow cell. Therein, the following PDE system has to be satisfied in the unit square domain $D = (0,1) \times (0,1)$
\begin{align*}
q(x,\omega) &= -a(x,\omega)\nabla v(x,\omega), \text{ for } x\in D,
\\ \nabla \cdot q(x,\omega) &= 0, \text{ for } x\in D ,
\end{align*}
for a.e. $\omega\in\Omega$, where $q$ is the Darcy velocity, $v$ is the hydrostatic pressure and $a$ is the permeability of the porous medium, which is modelled as a log-normal random field. More precisely, $\log(a)$ is a Gaussian random field with mean $\mu_Z=0$ and constant variance $\zeta_Z^2 = 1$. Moreover, $Z$ has an exponential type covariance function 
\begin{align*}
c(x,y) = \zeta_Z^2\exp\left(-\frac{\Vert x -y\Vert_1}{\lambda}\right),
\end{align*}
where $\lambda = 0.5$ denotes the correlation length. Again, the random field $Z$ is discretized by its KL expansion. The PDE system is coupled with the following boundary conditions
\begin{align}
\nu \cdot q(x,\omega) &= 0 \text{ for } x\in (0,1) \times \{0,1\},\label{2d bc 1}
\\ v(x,\omega) &= 1 \text{ for } x\in \{0\} \times (0,1),\label{2d bc 2}
\\ v(x,\omega) &= 0 \text{ for } x\in \{1\} \times (0,1),\label{2d bc 3}
\end{align}
for a.e. $\omega\in\Omega$, where $\nu$ denotes the derivative with respect to the normal direction on the boundary. Equation (\ref{2d bc 1}) impose that there is no flow across the horizontal boundaries and (\ref{2d bc 2}), (\ref{2d bc 3}) impose that there is inflow at the western boundary and outflow at the eastern boundary, respectively. The Darcy velocity $q$ is discretized by lowest order Raviart-Thomas mixed FEs, see \cite{Raviart77}. The pressure $v$ is discretized by piecewise constant elements. The grid is determined by the mesh size $h$ and consists of $2\cdot {1}/{h^2}$ uniform triangles. 
\\The failure event is based on the time that a particle requires to travel from the initial point $x_0 = (0, 0.5)^T$ to any other point on the boundary $\partial D$. Given the Darcy velocity $q_{h_{\ell}}(x,\omega)$, the particle path $x(t,\omega)$ has to satisfy the following ordinary differential equation
\begin{align*}
\frac{\partial}{\partial t}x(t,\omega) = q_{h_{\ell}}(x(t,\omega),\omega), \text{ } x(0,\omega) = x_0.
\end{align*} 
We approximate the particle path with the forward Euler discretization
\begin{align*}
x_{h_{\ell}}(t+\Delta t,\omega) = x(t,\omega) + \Delta t q_{h_{\ell}}(x(t,\omega),\omega),\text{ where } \Delta t = \frac{h_{\ell}}{2\Vert q_{h_{\ell}}(x(t,\omega),\omega)\Vert_2}.
\end{align*}
The travel time $\tau_{h_{\ell}}(\omega)\in[0,\infty)$ is defined as
\begin{align*}
\tau_{h_{\ell}}(\omega) = \underset{t>0}{\mathrm{argmin}}\text{ } x_{h_{\ell}}(t,\omega) \in \partial D.
\end{align*}
The approximation of the particle path is different to the procedure in \cite{Ullmann15} and, therefore, the estimated probability of failures differs slightly. Failure is defined as the event that $\tau_{h_{\ell}}$ is smaller than the threshold $\tau_0=0.03$. Hence, the respective LSF is defined as $G_{\ell}(\omega) := \tau_{h_{\ell}}(\omega)-\tau_0$. The reference solution of the probability of failure is $4.6730 \cdot 10^{-7}$ and is the estimated mean probability of failure over $100$ realizations of SuS with $N=10^4$ samples, mesh size $h={1}/{128}$, $\hat{p}_0=0.1$ and aCS as the MCMC method without burn-in. We note that SuS is a biased estimator and the relative bias scales as $\mathcal{O}(1/N)$ while the coefficient of variation scales as $\mathcal{O}(1/\sqrt{N})$ \cite{Au01}. The coefficient of variation of the $100$ probability of failure estimates is roughly $15\%$. Hence, we expect that the relative bias of the reference estimate is of order $10^{-2}$.
\\Figure \ref{flow_cell} shows a realization of a non-failure event and of a failure event. The figure displays the permeability $a$ and the respective solutions of the Darcy velocity $q_h$ for $h={1}/{128}$ and shows the particle paths which start at $x_0$ and their respective travel times. 
\\The probability of failure is estimated by SIS, MLSIS, SuS and MLSuS. For all methods, the estimation is performed for $N=250, 500, 1000$ samples and $N_s = 0.1\cdot N$ samples are considered for the small sample subset to decide if either bridging or tempering is performed in the update scheme of Section \ref{Update scheme}. For each parameter setting, the estimation is repeated $100$ times. For the multilevel methods, the sequence of mesh sizes is $h_\ell = 2^{-\ell-1}$ for $\ell=1,..,6$, i.e., the coarsest mesh size is $h_1={1}/{4}$ and the finest $h_6={1}/{128}$. The multi-level methods are applied with a level dependent dimension, where the parameter dimensions of the KL expansions are $n_1=10$, $n_2=20$, $n_3=40$, $n_4=80$ and $n_5=n_6=150$. SIS and MLSIS are performed for target coefficient of variations equal to $0.50$ and $1.00$. aCS and the vMFN distributions are considered as the MCMC methods without a burn-in. The parameter $c$ to define the number of seeds of the MCMC algorithm in (\ref{seed paramter}) is $c=0.1$. For SuS and MLSuS, aCS is considered as the MCMC method without a burn-in and the parameter $\hat{p}_0$ in (\ref{SuS parameter}) is $\hat{p}_0 =0.1$ or $\hat{p}_0 =0.25$. 

\subsubsection{Results}
\begin{figure}[htbp]
\centering
	\includegraphics[trim=0cm 0cm 0cm 0cm,scale=0.23]{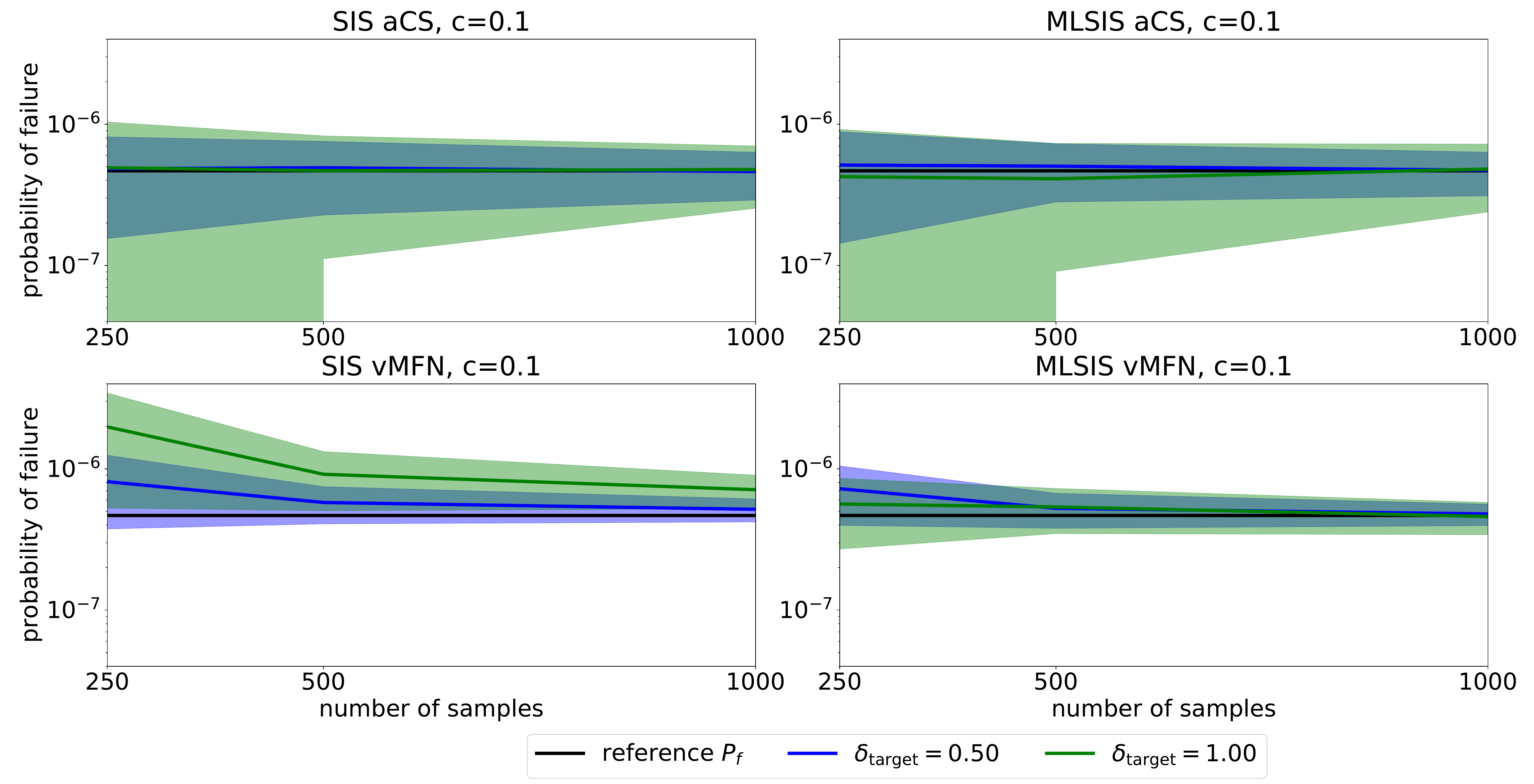}
		\caption[a]{Estimated probability of failure by SIS and MLSIS averaged over $100$ runs for $250, 500$ and $1000$ samples and $\delta_{\mathrm{target}}\in\{0.50, 1.00\}$. The coloured areas show the standard deviation of the estimates. The black lines show the reference estimate by Monte Carlo sampling. 1st row: aCS, $c=0.1$; 2nd row: vMFN, $c=0.1$; 1st column: SIS; 2nd column: MLSIS with level dependent dimension.}
		\label{2D-Diffusion: probability of failure + standard deviation}
\end{figure}

\begin{figure}[htbp]
\centering
	\includegraphics[trim=0cm 0cm 0cm 0cm,scale=0.23]{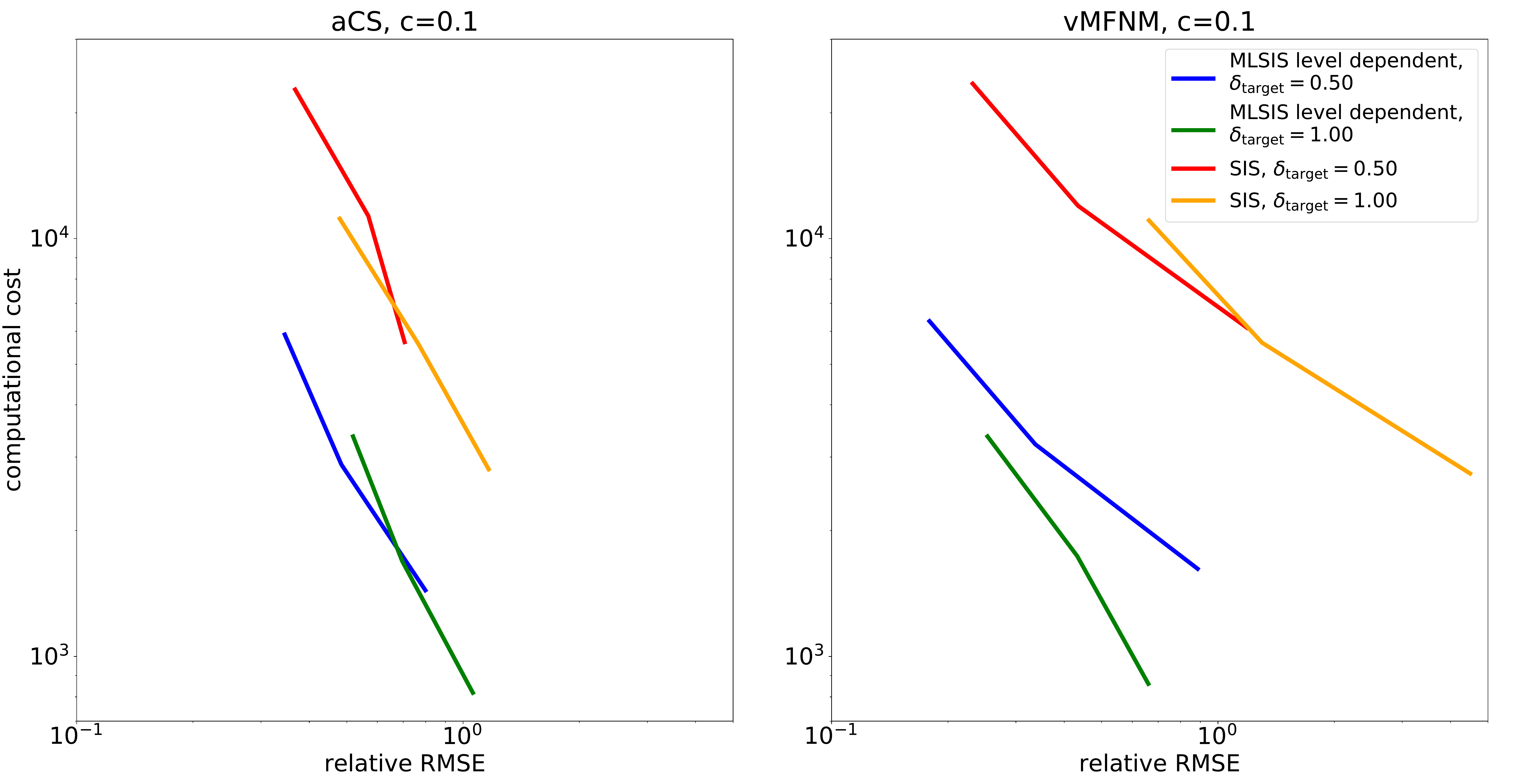}
		\caption[1D-Diffusion: RMSE and Costs]{Computational costs and relative RMSE of SIS and MLSIS averaged over $100$ runs for $250, 500$ and $1000$ samples and $\delta_{\mathrm{target}}\in\{0.50, 1.00\}$. 1st column: aCS, $c=0.1$; 2nd column: vMFN, $c=0.1$.}
		\label{2D-Diffusion: RMSE and Costs}
\end{figure}

\begin{figure}[h!]
	\centering
	\includegraphics[trim=0cm 0cm 0cm 0cm,scale=0.23]{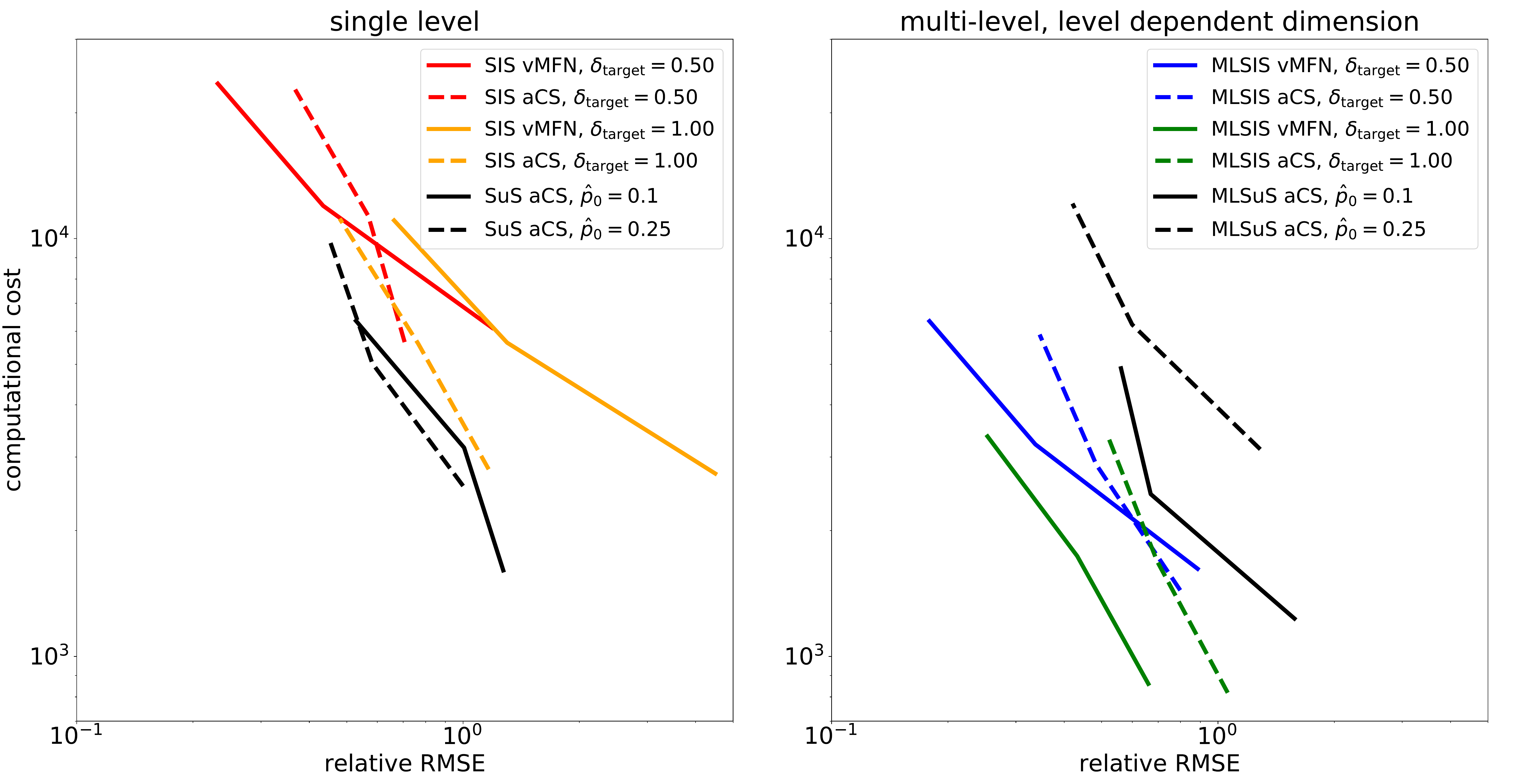}
		\caption[1D-Diffusion: Comparison RMSE and Costs]{Computational costs and relative RMSE of SIS, MLSIS, SuS and MLSuS averaged over $100$ runs for $250, 500$ and $1000$ samples. SIS and MLSIS are considered with aCS and vMFN, $c=0.1$ and $\delta_{\mathrm{target}}\in\{0.50, 1.00\}$. SuS and MLSuS are considered with aCS and $\hat{p}_0 =0.1$ or $\hat{p}_0 =0.25$. 1st column: single level; 2nd column: multi-level with level dependent dimension.}
		\label{2D-Diffusion: Comparison RMSE and Costs}
\end{figure}

Figure \ref{2D-Diffusion: probability of failure + standard deviation} shows the estimated mean probability of failure calculated by SIS and MLSIS plus/minus its standard deviation. The estimates of the means are in accordance with the reference solution and the bias and standard deviation decrease with an increasing number of samples. The standard deviation is smaller for a smaller target coefficient of variation. As for the 1D problem, we observe that applying the independent sampler with the vMFN distribution yields a smaller standard deviation than applying aCS.
\\Figure \ref{2D-Diffusion: RMSE and Costs} shows the relative RMSE on the horizontal axis and the computational costs on the vertical axis for the SIS and MLSIS estimator. The costs are calculated based on the formula given in (\ref{Cost levels}) for $L=6$ and $d=2$. Again, SIS and MLSIS yield the same range of the relative RMSE but the computational costs are lower for MLSIS. Considering the computational costs shown in Figure \ref{2D-Diffusion: RMSE and Costs}, we can save around $61\%$ of the computational costs if we apply MLSIS for the estimation. This is the same level of savings as in the 1D problem. However in the 2D problem, less level updates have to be performed as in the 1D problem setting. We expect that even more computational costs can be saved with MLSIS if we increase the highest discretization level.
\\Figure \ref{2D-Diffusion: Comparison RMSE and Costs} shows the relative RMSE and computational costs of SIS, MLSIS, SuS and MLSuS. We observe that SuS yields the same relative RMSE as SIS with aCS. However, SuS requires less computational costs. If we consider SIS with vMFN, the relative RMSE is smaller compared to SuS but the computational costs are higher for SIS. For the multilevel methods, we observe that MLSIS with $\delta_{\mathrm{target}}=1.00$ yields a smaller relative RMSE and requires less computational costs than MLSuS. This observation holds for both MCMC algorithms. In the 1D problem, we only observe that MLSIS with sampling from the vMFN distribution yields a more efficient estimator than MLSuS. For SuS and MLSuS, $\hat{p}_0 =0.25$ yields higher computational costs and a slightly smaller relative RMSE than $\hat{p}_0 =0.1$ since more intermediate failure domains are considered for $\hat{p}_0 = 0.25$. 

\section{Conclusion and Outlook}\label{Section Conclusion}
Motivated by the nestedness issue of Multilevel Subset Simulation, we implement Multilevel Sequential Importance Sampling to estimate the probability of rare events. We assume that the underlying limit state function depends on a discretization parameter $\ell$. MLSIS samples a sequence of non-zero density functions that are adaptively chosen such that each pair of subsequent densities are only slightly different. Therefore, nestedness is not an issue for MLSIS. We combine the smoothing approach of the indicator function in \cite{Papaioannou16} and the multilevel idea in \cite{Latz18}. This yields a two-fold adaptive algorithm which combines tempering and bridging sequences in a clever way to reduce computational costs. Moreover, we apply the level dependent dimension approach of \cite{Ullmann15} to reduce variances between consecutive accuracy levels of the limit state function. This leads to more tempering updates on coarse levels and reduces computational costs. Another contribution of our work is the von Mises Fisher Nakagami distribution as a proposal density in an independent Markov chain Monte Carlo algorithm. This leads to an efficient MCMC algorithm even in high dimensions. 
\\In numerical experiments in 1D and 2D space, we show for our experiments that MLSIS has a lower computational cost than SIS at any given error tolerance. For both experiments, MLSIS with the von Mises Fisher Nakagami distribution leads to lower computational cost than Multilevel Subset Simulation for the same accuracy. However, MLSIS with adaptive conditional sampling leads only for the 2D experiment to lower computational cost than Multilevel Subset Simulation for the same accuracy. The results also show that applying the von Mises Fisher Nakagami distribution as a proposal density in the MCMC algorithm reduces the bias and coefficient of variation of the MLSIS estimator compared to applying adaptive conditional sampling as the MCMC algorithm. 
\\The bridging approach can also handle more general assumptions on the approximation sequence of the limit state function. For instance, the approximation sequence can arise within a multi-fidelity setting. Therein, bridging is applied to transfer samples between a low fidelity model and a high fidelity model. 
\\Instead of using SIS to shift samples to the failure region, we plan to apply the Ensemble Kalman Filter for inverse problems as a particle based estimator for the probability of failure. In this case, the reliability problem is formulated as an inverse problem.

\bibliographystyle{plain}      
\bibliography{literatur}   

\begin{thebibliography}{10}

\bibitem{Abramowitz64}
M.~Abramowitz and I.~A. Stegun.
\newblock {\em Handbook of Mathematical Functions with Formulas, Graphs, and
  Mathematical Tables}.
\newblock U.S. Department of Commerce, National Bureau of Standards, 1964.

\bibitem{Agapiou17}
S.~Agapiou, O.~Papaspiliopoulos, D.~Sanz-Alonso, and A.~M. Stuart.
\newblock Importance sampling: Intrinsic dimension and computational cost.
\newblock {\em Statistical Science}, 32(3):405--431, 2017.

\bibitem{Agarwal17}
A.~Agarwal, S.~De~Marco, E.~Gobet, and G.~Liu.
\newblock Rare event simulation related to financial risks: efficient
  estimation and sensitivity analysis.
\newblock Working paper.

\bibitem{Au01}
S.-K. Au and J.~L. Beck.
\newblock Estimation of small failure probabilities in high dimensions by
  subset simulation.
\newblock {\em Probabilistic Engineering Mechanics}, 16(4):263--277, 2001.

\bibitem{Au14}
S.-K. Au and Y.~Wang.
\newblock {\em Engineering risk assessment with subset simulation}.
\newblock John Wiley \& Sons, 2014.

\bibitem{Beskos17}
A.~Beskos, A.~Jasra, K.~Law, R.~Tempone, and Y.~Zhou.
\newblock Multilevel sequential {Monte Carlo} samplers.
\newblock {\em Stochastic Processes and their Applications}, 127(5):1417--1440,
  2017.

\bibitem{Beskos13}
A.~Beskos, A.~Jasra, and A.~Thiery.
\newblock On the convergence of adaptive sequential {Monte Carlo} methods.
\newblock {\em The Annals of Applied Probability}, 26(2):1111--1146, 2013.

\bibitem{Botev12}
Z.~Botev and D.~Kroese.
\newblock Efficient {Monte Carlo} simulation via the generalized splitting
  method.
\newblock {\em Statistics and Computing}, 22(1):1--16, 2012.

\bibitem{Bouhlel15}
N.~Bouhlel and A.~Dziri.
\newblock Maximum likelihood parameter estimation of {Nakagami}-gamma shadowed
  fading channels.
\newblock {\em IEEE Communications Letters}, 19(4):685--688, 2015.

\bibitem{Braess07}
D.~Braess.
\newblock {\em Finite Elements: Theory, Fast Solvers, and Applications in Solid
  Mechanics}.
\newblock Cambridge University Press, 3 edition, 2007.

\bibitem{Chopin01}
N.~Chopin.
\newblock A sequential particle filter method for static models.
\newblock {\em Biometrika}, 89(3):539--551, 2002.

\bibitem{Cornaton08}
F.~J. Cornaton, Y.‐J. Park, S.~D. Normani, E.~A. Sudicky, and J.~F. Sykes.
\newblock Use of groundwater lifetime expectancy for the performance assessment
  of a deep geologic waste repository.
\newblock {\em Water Resources Research}, 44(4), 2008.

\bibitem{Cotter13}
S.~L. Cotter, G.~O. Roberts, A.~M. Stuart, and D.~White.
\newblock {MCMC} methods for functions: Modifying old algorithms to make them
  faster.
\newblock {\em Statistical Science}, 28(3):424--446, 2013.

\bibitem{Cerou12}
F.~Cérou, P.~Del~Moral, T.~Furon, and A.~Guyader.
\newblock Sequential {Monte Carlo} for rare event estimation.
\newblock {\em Statistics and Computing}, 22(3):795--808, 2012.

\bibitem{Dashti2017}
M.~Dashti and A.~M. Stuart.
\newblock {\em The {Bayesian} Approach to Inverse Problems}, pages 311--428.
\newblock Springer, 2017.

\bibitem{Angelis15}
M.~de~Angelis, E.~Patelli, and M.~Beer.
\newblock Advanced line sampling for efficient robust reliability analysis.
\newblock {\em Structural Safety}, 52:170--182, 2015.

\bibitem{Moral17}
P.~Del~Moral, A.~Jasra, K.~Law, and Y.~Zhou.
\newblock Multilevel sequential {Monte Carlo} samplers for normalizing
  constants.
\newblock {\em ACM Transactions on Modeling and Computer Simulation},
  27(3):20:1--20:22, 2017.

\bibitem{Kiureghian86}
A.~Der~Kiureghian and P.-L. Liu.
\newblock Structural reliability under incomplete probability information.
\newblock {\em Journal of Engineering Mechanics}, 112(1):85--104, 1986.

\bibitem{Doucet09}
A.~Doucet and A.~Johansen.
\newblock A tutorial on particle filtering and smoothing: Fifteen years later.
\newblock {\em Handbook of Nonlinear Filtering}, 12(3):656--704, 2009.

\bibitem{Elfverson16}
D.~Elfverson, F.~Hellman, and A.~M\r{a}lqvist.
\newblock A multilevel {Monte Carlo} method for computing failure
  probabilities.
\newblock {\em SIAM/ASA Journal on Uncertainty Quantification}, 4(1):312--330,
  2016.

\bibitem{Fishman96}
G.~S. Fishman.
\newblock {\em {Monte Carlo}: Concepts, Algorithms and Applications}.
\newblock Springer, New York, NY, 1996.

\bibitem{Gerber2019}
M.~Gerber, N.~Chopin, and N.~Whiteley.
\newblock Negative association, ordering and convergence of resampling methods.
\newblock {\em The Annals of Statistics}, 47(4):2236--2260, 2019.

\bibitem{Geyer19}
S.~Geyer, I.~Papaioannou, and D.~Straub.
\newblock Cross entropy-based importance sampling using {Gaussian} densities
  revisited.
\newblock {\em Structural Safety}, 76:15 -- 27, 2019.

\bibitem{Ghanem91}
R.~Ghanem and P.~Spanos.
\newblock {\em Stochastic Finite Elements: A Spectral Approach}.
\newblock Springer-Verlag, New York, 1991.

\bibitem{Gibbs02}
J.~W. Gibbs.
\newblock {\em Elementary Principles in Statistical Mechanics}.
\newblock Charles Scribner's Sons, 1902.

\bibitem{Giles15}
M.~B. Giles.
\newblock Multilevel {Monte Carlo} methods.
\newblock {\em Acta Numerica}, 24:259–328, 2015.

\bibitem{Glasserman99}
P.~Glasserman, P.~Heidelberger, P.~Shahabuddin, and T.~Zajic.
\newblock Multilevel splitting for estimating rare event probabilities.
\newblock {\em Operations Research}, 47(4):585--600, 1999.

\bibitem{Hastings70}
W.~K. Hastings.
\newblock {Monte Carlo} sampling methods using {Markov} chains and their
  applications.
\newblock {\em Biometrika}, 57(1):97--109, 1970.

\bibitem{Hohenbichler81}
M.~Hohenbichler and R.~Rackwitz.
\newblock Non-normal dependent vectors in structural safety.
\newblock {\em Journal of the Engineering Mechanics Division},
  107(6):1227--1238, 1981.

\bibitem{Jasra11}
A.~Jasra, D.~A. Stephens, A.~Doucet, and T.~Tsagaris.
\newblock Inference for {L\'{e}vy}-driven stochastic volatility models via
  adaptive sequential {Monte Carlo}.
\newblock {\em Scandinavian Journal of Statistics}, 38(3):1--22, 2011.

\bibitem{Kahn53}
H.~Kahn and A.~W. Marshall.
\newblock Methods of reducing sample size in {Monte Carlo} computations.
\newblock {\em Journal of the Operations Research Society of America},
  1(5):263--278, 1953.

\bibitem{Katafygiotis08}
L.~S. Katafygiotis and K.~M. Zuev.
\newblock Geometric insight into the challenges of solving high-dimensional
  reliability problems.
\newblock {\em Probabilistic Engineering Mechanics}, 23(2):208 -- 218, 2008.

\bibitem{Klenke2014}
A~Klenke.
\newblock {\em {Probability theory: A comprehensive course}}.
\newblock Springer, London, 2014.

\bibitem{Koutsourelakis09}
P.~S. Koutsourelakis.
\newblock A multi-resolution, non-parametric, {Bayesian} framework for
  identification of spatially-varying model parameters.
\newblock {\em Journal of Computational Physics}, 228(17):6184--6211, 2009.

\bibitem{Koutsourelakis04}
P.~S. Koutsourelakis, H.~J. Pradlwarter, and G.~I. Schuëller.
\newblock Reliability of structures in high dimensions, part i: algorithms and
  applications.
\newblock {\em Probabilistic Engineering Mechanics}, 19(4):409 -- 417, 2004.

\bibitem{Lacaze15}
S.~Lacaze, L.~Brevault, S.~Missoum, and M.~Balesdent.
\newblock Probability of failure sensitivity with respect to decision
  variables.
\newblock {\em Structural and Multidisciplinary Optimization}, 52(2):375--381,
  2015.

\bibitem{Latz18}
J.~Latz, I.~Papaioannou, and E.~Ullmann.
\newblock Multilevel $\mathrm{Sequential}^2$ {Monte Carlo} for {Bayesian}
  inverse problems.
\newblock {\em Journal of Computational Physics}, 368:154--178, 2018.

\bibitem{McLachlan05}
G.~McLachlan and D.~Peel.
\newblock {\em {ML} Fitting of Mixture Models}, chapter~2, pages 40--80.
\newblock John Wiley and Sons, Ltd, 2005.

\bibitem{Melchers18}
R.~E. Melchers and A.~T. Beck.
\newblock {\em Structural Reliability Analysis and Prediction}.
\newblock John Wiley and Sons, 3 edition, 2018.

\bibitem{Metropolis53}
H.~Metropolis, A.~Rosenbluth, M.~Rosenbluth, A.~Teller, and E.~Teller.
\newblock Equations of state calculations by fast computing machine.
\newblock {\em The Journal of Chemical Physics}, 21(6):1087--1092, 1953.

\bibitem{Morio15}
J.~Morio and M.~Balesdent.
\newblock {\em Estimation of Rare Event Probabilities in Complex Aerospace and
  other Systems}.
\newblock Woodhead Publishing, 1 edition, 2015.

\bibitem{Nakagami60}
M.~Nakagami.
\newblock The m-distribution, a general formula of intensity distribution of
  rapid fading.
\newblock In W.C. Hoffman, editor, {\em Statistical Methods in Radio Wave
  Propagation}, pages 3--36. Pergamon, 1960.

\bibitem{Noseck08}
U.~Noseck, D.~Becker, C.~Fahrenholz, E.~Fein, J.~Fl{\"u}gge, K.-P. Kr{\"o}hn,
  J.~M{\"o}nig, I.~M{\"u}ller-Lyda, T.~Rothfuchs, A.~R{\"u}bel, and J.~Wolf.
\newblock {\em Assessment of the long-term safety of repositories}.
\newblock Gesellschaft für Anlage und Reaktorsicherheit (GRS) mbH, 2008.

\bibitem{Papaioannou15}
I.~Papaioannou, W.~Betz, K.~Zwirglmaier, and D.~Straub.
\newblock {MCMC} algorithms for subset simulation.
\newblock {\em Probabilistic Engineering Mechanics}, 41:89--103, 2015.

\bibitem{Papaioannou18}
I.~Papaioannou, K.~Breitung, and D.~Straub.
\newblock Reliability sensitivity estimation with sequential importance
  sampling.
\newblock {\em Structural Safety}, 75:24 -- 34, 2018.

\bibitem{Papaioannou19}
I.~Papaioannou, S.~Geyer, and D.~Straub.
\newblock Improved cross entropy-based importance sampling with a flexible
  mixture model.
\newblock {\em Reliability Engineering \& System Safety}, 191:106564, 2019.

\bibitem{Papaioannou16}
I.~Papaioannou, C.~Papadimitriou, and D.~Straub.
\newblock Sequential importance sampling for structural reliability analysis.
\newblock {\em Structural Safety}, 62:66--75, 2016.

\bibitem{Peherstorfer18}
B.~Peherstorfer, B.~Kramer, and K.~Willcox.
\newblock Multifidelity preconditioning of the cross-entropy method for rare
  event simulation and failure probability estimation.
\newblock {\em SIAM/ASA Journal on Uncertainty Quantification}, 6(2):737--761,
  2018.

\bibitem{Rackwitz01}
R.~Rackwitz.
\newblock Reliability analysis—a review and some perspectives.
\newblock {\em Structural Safety}, 23(4):365--395, 2001.

\bibitem{Raviart77}
P.~A. Raviart and J.~M. Thomas.
\newblock {\em A Mixed Finite Element Method for Second Order Elliptic
  Problems}, volume 606, pages 292--315.
\newblock Springer-Verlag, Berlin, New York, 1977.

\bibitem{roberts2001}
G.~O. Roberts and J.~S. Rosenthal.
\newblock Optimal scaling for various {Metropolis-Hastings} algorithms.
\newblock {\em Statistical Science}, 16(4):351--367, 2001.

\bibitem{Rubinstein16}
R.~Y. Rubinstein and D.~P. Kroese.
\newblock {\em Simulation and the {Monte Carlo} Method: Concepts, Algorithms
  and Applications}.
\newblock John Wiley and Sons, 3 edition, 2016.

\bibitem{Ullmann15}
E.~Ullmann and I.~Papaioannou.
\newblock Multilevel estimation of rare events.
\newblock {\em SIAM/ASA Journal on Uncertainty Quantification}, 3(1):922--953,
  2015.

\bibitem{Wang19}
Z.~Wang, M.~Broccardo, and J.~Song.
\newblock Hamiltonian {Monte Carlo} methods for subset simulation in
  reliability analysis.
\newblock {\em Structural Safety}, 76:51--67, 2019.

\bibitem{Wang16}
Z.~Wang and J.~Song.
\newblock Cross-entropy-based adaptive importance sampling using von
  {Mises-Fisher} mixture for high dimensional reliability analysis.
\newblock {\em Structural Safety}, 59:42--52, 2016.

\end{thebibliography}

\end{document}